\documentclass[12pt]{article}
\usepackage{mathrsfs}
\usepackage{amsmath}
\usepackage{amsthm}
\usepackage[all]{xy}
\usepackage{titlesec}
\usepackage{amssymb}
\newtheorem{theorem}{Theorem}[section]
\newtheorem{definition}[theorem]{Definition}
\newtheorem{lemma}[theorem]{Lemma}
\newtheorem{proposition}[theorem]{Proposition}

\newtheorem{Mthm}{Main Theorem}

\theoremstyle{plain}

\newcommand{\eb}{\begin{enumerate}}
\newcommand{\ee}{\end{enumerate}}
\newcommand{\ebx}{\begin{equation*}}
\newcommand{\eex}{\end{equation*}}

\newcommand{\beq}{\begin{eqnarray*}}
\newcommand{\eeq}{\end{eqnarray*}}

\makeatletter

\newcommand{\Rmnum}[1]{\expandafter\@slowromancap\romannumeral #1@}
\makeatother

\begin{document}
\setcounter{page}{1}
\title{Admissible pairs of Hermitian symmetric spaces in the perspective of the theory of varieties of minimal rational tangents
}


\author{Yunxin Zhang}




\date{}
\maketitle

\begin{abstract}
We study a pair $(\mathcal{S}_0,\mathcal{S})$ of irreducible Hermitian Symmetric Spaces of compact type (cHSS) in this paper, with the first aim being classifying all the \textsl{admissible pairs} $(\mathcal{S}_0,\mathcal{S})$. This notion is a natural generalization of \textsl{the pairs of} \textsl{sub-diagram type} originated by Jaehyun Hong and Ngaiming Mok (\cite{[HoM 10]}). Based on this classification, we partially solve the rigidity problem for the admissible pairs $(\mathcal{S}_0,\mathcal{S})$ which was raised by Mok and Zhang (2014) (\cite{[MoZ 14]}), culminating in determining a sufficient condition for the pairs being non-rigid and proving that \textsl{special pairs}, which show up in the classification procedure, are algebraic, as a weaker result than being rigid. However, whether special pairs are rigid or not remains unknown and needs further investigation in the framework of VMRT theory.
\end{abstract}



\section{Introduction and Motivation}
\label{intro}

It is well known that a Rational Homogeneous Space (RHS) $X=G/P$ is a projective manifold.  For $X$ of Picard number 1,  $\Gamma(X, \mathcal{O}(1)): X\hookrightarrow\mathbb{P}^N$  realises the so called first canonical embedding where $\mathcal{O}(1)$ denotes the generator of the 1-dimensional Picard group.  Consider the preimages $C$ of the projective lines $L\subset\mathbb{P}^N$  under the first canonical embedding, we call $C$ the \textsl{minimal rational curve} (MRC). The collection of tangent vectors $[T_xC]$ of all the minimal rational curves $C$ passing through an arbitrary fixed point $x\in X$  is called the variety of minimal rational tangents (VMRT) of $X$, denoted by $\mathscr{C}_x(X)$. It's a projective submanifold $\mathscr{C}_x(X)\subset\mathbb{P}(T_xX)$, we may also refer to the associated affine submanifold $\widetilde{\mathscr{C}}_x(X)\subset T_xX$. Note that since $X$ is homogeneous,  then $\mathscr{C}_x(X)$ is projectively equivalent to $\mathscr{C}_y(X)$ for any $x,y\in X$. For a comprehensive survey of VMRT theory, the reader may refer to \cite{[Mok 08]} or \cite{[HwM 99(b)]}. In a series of joint works, Jun-Muk Hwang and Ngaiming Mok have developed a geometric theory of VMRT on uniruled manifold, which covers rational homogeneous spaces. This theory turns out to be very effective in solving some classical algebro-geometric problems, including rigidity of rational homogeneous spaces under K$\ddot{\textnormal{a}}$hler deformation, Larzasfeld problem etc. (cf. \cite{[HwM 98]}, \cite{[HwM 99(a)]}, \cite{[HwM 02]}, \cite{[HwM 05]}). In addition, the joint work of Jaehyun Hong and Ngaiming Mok provided another aspect of the application of the VMRT theory, in the sense of characterising  standard embedding between rational homogeneous spaces of certain type,  the precise statement of their main result is as follows:

\begin{theorem}
\label{Thm 1.1}
\textnormal{(cf. theorem 1.2 in [HoM 10])} Let $X_0=G_0/P_0,X=G/P$ be irreducible rational homogeneous spaces (RHS) associated to a long simple root determined by marked Dynkin diagrams $(\mathcal{D}(G_0),\gamma_0),(\mathcal{D}(G),\gamma)$ respectively. Suppose  $\mathcal{D}(G_0)$ is obtained from a sub-diagram of $\mathcal{D}(G)$ with $\gamma_0$ being identified with $\gamma$. If $X_0$ is non-linear and $f:U\rightarrow X$ is a holomorphic embedding from a connected open subset $U\subset X_0$ into $X$ which respects VMRT at a general point $x\in U$, then $f$ is the restriction of a standard embedding of $X_0$ into $X$.
\end{theorem}

The pairs $(X_0,X)$ of RHS as described in the above theorem is said to be of \textsl{sub-diagram type}. The choice of a sub-diagram $\mathcal{D}(G_0)$ of $\mathcal{D}(G)$ naturally induces an equivariant embedding $i:X_0=G_0/P_0\hookrightarrow X=G/P$. For any $g\in G$, we say the composition $g\circ i:X_0\hookrightarrow X$ is a \textsl{standard embedding}, which satisfies some crucial properties in the perspective of VMRT theory such as \textsl{VMRT respecting} and \textsl{preservation of minimal rational curves}. In this sense, we say pair $(X_0,X)$ of sub-diagram type is an admissible pair, whose formal definition is as below\footnote{The following three definitions \ref{Def 1.2 for paper 2}, \ref{Def 1.4}, \ref{Def 1.4'}  are originated in \cite{[Zhang 14]}}:

\begin{definition}
\label{Def 1.2 for paper 2}
A pair of rational homogeneous spaces $(G_0/P_0,G/P)$ of Picard number 1 associated to simple roots $(\gamma_0,\gamma)$ is said to be an admissible pair if and only if there exists an equivariant holomorphic embedding $i:X_0\hookrightarrow X$ such that $i_*:H_2(X_0,\mathbb{Z})\stackrel{\cong}{\rightarrow}H_2(X,\mathbb{Z})$ and such that $i$ respects VMRTs, i.e., d$i(\mathscr{C}_x(X_0))=$d$i(\mathbb{P}T_x(X_0))\cap\mathscr{C}_{i(x)}(X)$ for every $x\in X_0$. Such $i$ is also said to be a standard embedding; the image $i(X_0)\subset X$ is said to be a standard model.
\end{definition} 

\textbf{Remark} Admissible pairs are transitive, i.e., if $(X'_0,X_0),(X_0,X)$ are both admissible, then so is $(X'_0,X)$. 

Following the notion of admissible pairs, the definition of \textsl{sub-VMRT strucure} naturally arises as:

\begin{definition}
\label{Def 1.4}
Let $(X_0,X)$ be an admissible pair of rational homogeneous manifolds of Picard number 1, $W\subset X$ be an open subset and $S\subset W$ be a complex submanifold. Consider the fibered space $\pi:\mathscr{C}(X)\rightarrow X$ of varieties of minimal rational tangents on $X$. Define $\mathscr{C}(S)\subset\mathscr{C}(X)|_S$ by $\mathscr{C}_x(S):=\mathscr{C}_x(X)\cap\mathbb{P}T_x(S)$, we say that $S\subset U$ inherits a sub-VMRT structure modeled on $(X_0,X)$ if and only if for every point $x\in S$ there exists a neighborhood $U$ of $x\in S$ and a trivialization of the holomorphic projective bundle $\mathbb{P}T(X)|_U$ given by $\Phi:\mathbb{P}T(X)|_U\stackrel{\cong}{\rightarrow}\mathbb{P}T_o(X)\times U$ such that (1) $\Phi(\mathscr{C}(X)|_U)=\mathscr{C}_o(X)\times U$ and (2) $\Phi(\mathscr{C}(S)|_U)=\mathscr{C}_o(X_0)\times U$.
\end{definition}

Note that the definition that $S\subset W$ inherits a sub-VMRT structure modeled on the admissible pair $(X_0,X)$ can be reformulated as requiring:

$(\dag)$ For any $x\in S$, there exists a projective linear isomorphism $\Lambda_x:\mathbb{P}T_x(X)\stackrel{\cong}{\rightarrow}\mathbb{P}T_o(X)$ such that $\Lambda_x(\mathscr{C}_x(X))=\mathscr{C}_o(X)$ and $\Lambda_x(\mathscr{C}_x(S))=\mathscr{C}_o(X_0)$.

This is the case because given $(\dag)$, the set of fiberwise projective linear isomorphisms $\Lambda_x:\mathbb{P}T_x(X)\stackrel{\cong}{\rightarrow}\mathbb{P}T_o(X)$ satisfying the two requirements in $(\dag)$ forms a holomorphic fiber bundle over $S$, hence admitting at each $x\in S$ a holomorphic section on some neighborhood $U$ of $x\in S$.

The above definition provides us with a notion of imposing on submanifold a standard sub-VMRT structure, it is natural to ask if a standard sub-VMRT structure is sufficient for the submanifold itself to be some standard model $X_0$. This question inspires the following notion of \textsl{rigidity}:

\begin{definition}
\label{Def 1.4'}
An admissible pair $(X_0,X)$ of RHS of Picard number 1 is said to be rigid if and only if any locally closed complex submanifold $S\subset X$ inheriting a sub-VMRT structure modeled on $(X_0,X)$ must necessarily be an open subset of some standard model in $X$. In this case we also say that sub-VMRT structures modeled on $(X_0,X)$ are rigid.
\end{definition}

Then problem arises as whether we can give a sufficient and necessary condition for the admissible pairs $(X_0,X)$ to be rigid. The joint work of Hong and Mok (2010) as summarized in Theorem \ref{Thm 1.1} serves as the first achievement for this purpose. However, it turns out to be a huge challenge to even determine all the admissible pairs $(X_0,X)$ of RHS, let alone investigating the rigidity. Towards solving this problem, it is natural to first of all restrict our consideration to the scope of cHSS. The first task comes as classifying all the admissible pairs $(\mathcal{S}_0,\mathcal{S})$ of cHSS:

\begin{Mthm}
\label{Mthm 1 for paper 2}
The admissible pair $(\mathcal{S}_0,\mathcal{S})$ of cHSS where $\mathcal{S}_0$ is non-linear falls into at least one of following four categories: 
\begin{enumerate}
\item \textsl{pairs of sub-diagram type}

\item 
\begin{enumerate}
\item pairs of deletion type
\item $(Q^n,Q^m),\ \ n-m\equiv 1$(mod $2$);
\end{enumerate}

\item \textsl{special pairs}: 
\begin{enumerate}
\item$(G^{\Rmnum{3}}(n,n),G(r,s)),\ \ (3\leq n\leq $min$\{r,s\}$);
\item $(G(r,s),G^{\Rmnum{2}}(n,n))$, ($r,s\geq 3,\ \ r+s\leq n$);
\item $ (G(4,2),\Rmnum{5})$; 
\item $(G(5,2),\Rmnum{6})$;
\item $(G(6,2),\Rmnum{6})$; 
\item $(G(3,3),\Rmnum{6})$;
\item $(G^{\Rmnum{2}}(6,6),\Rmnum{6})$;
where $\Rmnum{5},\Rmnum{6}$\ \ denote $E_6/(\textnormal{Spin}(10)\times U(1)),E_7/(E_6\times U(1))$,\ \ respectively.
\end{enumerate}

\item By virtue of transitivity of admissible pair (see remark after Def \ref{Def 1.2 for paper 2}), there exists some cHSS $\mathcal{S}'$, such that both $(\mathcal{S}_0,\mathcal{S}')$ and $(\mathcal{S}',\mathcal{S})$ are admissible pairs which fall into one of the above three categories.
\end{enumerate}
\end{Mthm}

There appears the notion of \textsl{pair of deletion type} in category 2 in the above theorem, hyperquadrics $(Q^n,Q^m)$ of even dimensional difference, i.e., $m-n$ being an even number, provide typical examples for such kind of pairs. The formal definition of deletion type pairs comes as follows:

\begin{definition}
\label{Def 1.3 for paper 2}
Let the cHSSs $\mathcal{S}_0,\mathcal{S}',\mathcal{S}$ be associated to $(\mathcal{D}_0,\gamma_0),(\mathcal{D}',\gamma'),(\mathcal{D},\gamma)$ respectively. Suppose that $(\mathcal{S}_0,\mathcal{S}')$ is either of sub-diagram type (cf. Thm \ref{Thm 1.1}) or identical, and $(\mathcal{D}',\gamma')$ is obtained from $(\mathcal{D},\gamma)$ by deleting a chain whose one end node is attached to $\gamma'$ and other end node is $\gamma$, then we say the pair $(\mathcal{S}_0,\mathcal{S})$ is of deletion type. Here a chain means the Dynkin diagram of $A$ type.  
\end{definition}

Note that we rule out the pairs $(\mathcal{S}_0,\mathcal{S})$ in which $\mathcal{S}_0$ is a projective space. This is well justified by an early joint work of Hong and Choe, which competely determined the pairs $(Z_{\textnormal{max}},\mathcal{S})$ where $Z_{\textnormal{max}}$ denotes the \textsl{maximal linear space} (cf. \cite{[HoC 04]} for definition). Their theorem says:

\begin{theorem}
\label{Thm 1.4 for paper 2}
Let $(Z_{\textnormal{max}},\mathcal{S})$ be admissible pair of cHSS where $Z_{\textnormal{max}}$ is a maximal linear space. Denote by $d_{\textnormal{max}}$ the dimension of $Z_{\textnormal{max}}$, then the dependence of $d_{\textnormal{max}}$ on $\mathcal{S}$ is shown in the table below:

\begin{table}
\caption{Hermitian Symmetric Spaces and the dimension of maximal linear subspaces cf. [HoC 04]}
\label{tab:1}
	\begin{tabular}{|l|l|l|l|l|l|l|l|}
	\hline
	$\mathcal{S}$ & $G(p,q)$ & $G^{\Rmnum{2}}(n,n)$ & $G^{\Rmnum{3}}(n,n)$ & $Q^{2k}$ & $Q^{2k+1}$ & $E_6$ & $E_7$\\
	\hline
	$d_{\textnormal{max}}$ & $p$ or $q$ & $3$ or $n-1$ & $1$ & $k$ & $k$ & $4$ or $5$ & $5$ or $6$ \\
	\hline
	\end{tabular}
\end{table}
	
\end{theorem}

It's not hard to see that the combination of our Main Theorem \ref{Mthm 1 for paper 2} and Thm \ref{Thm 1.4 for paper 2} completely classify all the admissible pairs $(\mathcal{S}_0,\mathcal{S})$ of cHSS. Based on this classification, up until now, we have a partial solution of rigidity problem in the scope of cHSS, in the sense that we prove a sufficient condition for those pairs which are non-rigid as contrast to Thm \ref{Thm 1.1} which is a sufficient condition for pairs to be rigid.

\begin{Mthm}
\label{Mthm 2 for paper 2}
An admissible pair $(\mathcal{S}_0,\mathcal{S})$ of cHSS is non-rigid if it is degenerate.
\end{Mthm}

The notion of \textsl{degenerate} or \textsl{non-degenerate} pairs is originated in [HoM 10], whose formal definition is as follows:


\begin{definition}
\label{Def 1.5}
\textnormal{(cf. [HoM 10])} Let $o\in \mathcal{S}_0$ be some fixed reference point and  $\mathbb{C}\alpha\in di(\widetilde{\mathscr{C}}_o(\mathcal{S}_0))=\widetilde{\mathscr{C}}_{i(o)}(\mathcal{S})\cap di(T_o(\mathcal{S}_0))$. Denote by $\sigma'$ the second fundamental form with respect to the flat connection in $T_{i(x)}(\mathcal{S})$ as a Euclidean space and define

\begin{center}
Ker$'(\sigma'):=\{u\in T_{\alpha}(\widetilde{\mathscr{C}}_{i(o)}(\mathcal{S})):\quad\sigma'(u,w)=0,\ \ \forall\ \ w\in T_{\alpha}(di(\widetilde{\mathscr{C}}_o(\mathcal{S}_0)))=di(T_o(\mathcal{S}_0))\cap T_{\alpha}(\widetilde{\mathscr{C}}_{i(o)}(\mathcal{S}))\}$
\end{center}
If Ker $(\sigma')$=$\mathbb{C}\alpha$ when $\mathcal{S}_0$ is not projective space, or Ker$(\sigma')\subset T_{\alpha}(di(\widetilde{\mathscr{C}}_o(\mathcal{S}_0)))$ when $\mathcal{S}_0$ is projective space, then we say that the pair $(\mathcal{S}_0,\mathcal{S})$ is non-degenerate. Otherwise we say it is degenerate.
\end{definition}

Note that for two admissible pairs of cHSS $(\mathcal{S}'_0,\mathcal{S}_0),(\mathcal{S}_0,\mathcal{S})$, if either of them is degenerate, so is the admissible pair $(\mathcal{S}'_0,\mathcal{S})$. As a direct result of Thm \ref{Thm 1.1}, the pairs in the first category in Main Theorem 1 is rigid. It can be checked that pairs of deletion type is degenerate, so is non-rigid as a result of Main Theorem 2. While the non-rigidity of $(Q^n,Q^m)$ with $m-n$ being odd number is immediate from considering any $n$-dimensional submanifold $S\subset Q^m$ which must inherit a sub-VMRT structure modeled on $(Q^n, Q^m)$, it remains for us to investigate the rigidity for the \textsl{special pairs} in the third category. However, we can only prove  the weaker result ``algebraic'' instead of rigid for special pairs (full details to be shown in chapter 2). Not until we can confirm the special pairs to be rigid or non-rigid that we completely solve the rigidity problem in the scope of cHSS, it needs further investigation from the perspective of VMRT theory in the future. 

Before concluding this introduction, we provide some results which are probably existing in literatures but not readily available. The reader may refer to \cite{[Mok 89]},\cite{[HwM 02]},\cite{[HwM 99(b)]} at some points. 

A cHSS $\mathcal{S}=G/P$ is homogeneous where $G=$Aut$(\mathcal{S})$ and $P$ is the maximal parabolic subgroup fixing some reference point $o\in\mathcal{S}$. Every cHSS is associated with a unique \textsl{Harish-Chandra decomposition} $\frak{g}=\frak{m}^-\oplus\frak{k}^{\mathbb{C}}\oplus\frak{m}^+$, satisfying $[\frak{m}^-,\frak{k}^{\mathbb{C}}]\subset\frak{m}^-$, $[\frak{m}^+,\frak{k}^{\mathbb{C}}]\subset\frak{m}^+$, $[\frak{k}^{\mathbb{C}},\frak{k}^{\mathbb{C}}]\subset\frak{k}^{\mathbb{C}}$, $[\frak{m}^+,\frak{m}^-]\subset\frak{k}^{\mathbb{C}}$, $\frak{m}^+$ is an Abelian subalgebra such that $\frak{m}^+\cong T_o(\mathcal{S})$, i.e., the holomorphic tangent space and $\frak{m}^-=\overline{\frak{m}^+}$, the bar $\overline{\cdot}$ denotes the complex conjugation. $\frak{g}=$Lie$(G)$ is the simple complex Lie algebra and $\frak{p}=\frak{m}^-\oplus\frak{k}^{\mathbb{C}}$=Lie$(P)$. The subalgebra $\frak{k}^{\mathbb{C}}$ has a $1$-dimensional center giving rise to an integrable complex structure on $\mathcal{S}$, and the corresponding reductive complex Lie group $K^{\mathbb{C}}$ acts as automorphism of the VMRT $\mathscr{C}_o(\mathcal{S})$ via the isotropy representation on $\frak{m}^+\cong T_o(\mathcal{S})$, we collect some facts concerning the structure of $\mathscr{C}_o(\mathcal{S})$ as follows:

\begin{proposition}[cf. \cite{[Mok 89]}]
\label{Prop 1.8 for paper 2}
Denote by $\frak{k}_s=[\frak{k}^{\mathbb{C}},\frak{k}^{\mathbb{C}}]$ the semisimple part, and $K_s$ the corresponding Lie group, then the VMRT $\mathscr{C}_o(\mathcal{S})=\displaystyle\bigcup_{k\in K_s}$(Ad$(k))([\alpha])$ where $\alpha\in T_o(\mathcal{S})$ is some highest weight vector w.r.t. the irreducible isotropy representation of $K_s$ on $T_o(\mathcal{S})\cong\frak{m}^+$ and Ad$(k)$ denotes the adjoint action. The tangent space of VMRT at the point $[\alpha]\in\mathbb{P}(T_o(\mathcal{S}))$ is $T_{[\alpha]}(\mathscr{C}_o(\mathcal{S}))=([\frak{k}_s,\alpha]+\mathbb{C}\alpha)/\mathbb{C}\alpha$. It gives rise to a direct decomposition as $T_o(\mathcal{S})=\mathbb{C}\alpha\oplus\mathscr{H}_{\alpha}\oplus\mathscr{N}_{\alpha}$ where $\mathscr{H}_{\alpha}=T_{[\alpha]}(\mathscr{C}_o(\mathcal{S}))$ and $\mathscr{N}_{\alpha}$ is the normal complement subspace, with the holomorphic sectional or bisectional curvature $R^h(\alpha,\overline{\alpha},\alpha,\overline{\alpha})=2$, $R^h(\alpha,\overline{\alpha},\xi,\overline{\xi})=1$ if $\xi\in\mathscr{H}_{\alpha}$ and $0$ if $\xi\in\mathscr{N}_{\alpha}$, $h$ is some normalised canonical metric on $\mathcal{S}$.
\end{proposition}

Consider the fiber bundle $\mathscr{C}(\mathcal{S}):=\displaystyle\bigcup_{x\in\mathcal{S}}\mathscr{C}_x(\mathcal{S})$, it is remarkable that this fiber bundle is flat, i.e., there exists a coordinate covering $U_{\alpha}$ such that the fiber bundle is a product on coordinate charts, i.e., $\mathscr{C}(\mathcal{S})|_{U_{\alpha}}=\mathscr{C}_o(\mathcal{S})\times U_{\alpha}$. Such coordinate chart is \textsl{Harish-Chandra coordinate} as a result of \textsl{Harish-Chandra embedding}. Since the vector $\alpha\in\widetilde{\mathscr{C}}_o(\mathcal{S})$ is the tangent vector to a unique minimal rational curve $L$ passing through $o$, i.e., $T_oL=\mathbb{C}\alpha$, the flatness of $\mathscr{C}(\mathcal{S})$ implies minimal rational curves are expressed as affine lines in terms of coordinate. We summarize the results as follows;

\begin{proposition}[cf. [Mok 89]]
\label{Prop 2' for paper 2}
$\mathcal{S}=G/P$ is a cHSS of dimension $m$, let $M^-,K^{\mathbb{C}},M^+$ denote the respective Lie groups of $\frak{m}^-,\frak{k}^{\mathbb{C}},\frak{m}^+$. The holomorphic mapping $F:M^-\times K^{\mathbb{C}}\times M^+\rightarrow G$ defined by $F(m^-,k,m^+)=m^-km^+$ is biholomorphism of $M^-\times K^{\mathbb{C}}\times M^+$ onto a dense open subset of $G$. In particular, the mapping $\eta:\frak{m}^+\hookrightarrow\mathcal{S}$ given by $\eta(\frak{m}^+)=$exp$(\frak{m}^+)P$ is a biholomorphism onto a dense open subset $U\subset\mathcal{S}$, $U\cong\mathbb{C}^{m}$ is a called a Harish-Chandra coordinate chart. The fiber bundle $\mathscr{C}(\mathcal{S}):=\displaystyle\bigcup_{x\in\mathcal{S}}\mathscr{C}_x(\mathcal{S})$ on $\mathcal{S}$ is flat in the sense that $\mathscr{C}(\mathcal{S})|_{U}=\mathscr{C}_o(\mathcal{S})\times U$. As a result, minimal rational curves are affine lines in terms of coordinates in $U$. 
\end{proposition}

The VMRT $\mathscr{C}_o(\mathcal{S})\subset\mathbb{P}(T_o\mathcal{S})$ as a projective manifold, called \textsl{embedded VMRT}, is obtained through an embedding $\phi:\mathcal{A}(\mathcal{S})\hookrightarrow\mathbb{P}^N$ where $\mathcal{A}(\mathcal{S})$ is called the \textsl{abstract VMRT of $\mathcal{S}$}. It turns out $\mathcal{A}(\mathcal{S})$ is itself a cHSS of rank $\leq 2$ (cf. [Mok 89]). The embedding $\phi:\mathcal{A}(\mathcal{S})\cong\mathscr{C}_o(\mathcal{S})\subset\mathbb{P}(T_o\mathcal{S})$ is realised by $\mathcal{O}(1)$ when $\mathcal{S}$ is other than Lagrangian Grassmannian $G^{\Rmnum{3}}(n,n)$, i.e., $\phi$ is first canonical embedding which sends minimal rational curves of $\mathcal{A}(\mathcal{S})$ to projective lines in $\mathbb{P}(T_o\mathcal{S})$. In the case of $G^{\Rmnum{3}}(n,n)$ whose abstract VMRT is nothing but $\mathbb{P}^{n-1}$, $\phi$ is realised by $\mathcal{O}(2)$, i.e., the Veronese embedding of $\mathbb{P}^{n-1}$ (cf. \cite{[Mok 99]}). The following lemma is just a direct corollary of a vector valued cubic polynomial which can be found in section 4.2 of \cite{[HwM 99(b)]}.

\begin{lemma}
\label{Lemma 1 of paper 2}
Suppose an embedded VMRT $\mathscr{C}_o(\mathcal{S})\subset\mathbb{P}(T_o\mathcal{S})$, then $\mathscr{C}_o(\mathcal{S})$ is the Zariki closure of a vector valued quadratic polynomial $\zeta$ defined on the entire tangent space $\mathscr{H}_{\alpha}=T_{[\alpha]}(\mathscr{C}_o(\mathcal{S}))\cong T_{\alpha}(\widetilde{\mathscr{C}}_o(\mathcal{S}))/\mathbb{C}\alpha$. Precisely, 
\begin{center}
$\zeta:T_{\alpha}(\widetilde{\mathscr{C}}_o(\mathcal{S}))\rightarrow T_{\alpha}(\widetilde{\mathscr{C}}_o(\mathcal{S}))\oplus \mathscr{N}_{\alpha}=T_o(\mathcal{S})$\\
\vspace{4mm}
$\zeta(\xi)=\alpha+\xi+\sigma'(\xi,\xi),\ \ \forall\xi\in \mathscr{H}_{\alpha}$
\end{center}where $\sigma':  \textnormal{Sym}^2(T_{\alpha}(\widetilde{\mathscr{C}}_o(\mathcal{S})))\rightarrow \mathscr{N}_{\alpha}$ is the second fundamental form of $\widetilde{\mathscr{C}}_o(\mathcal{S})\backslash\{0\}$ in $T_o(\mathcal{S})$ w.r.t. Euclidean flat metric. 

\end{lemma}

Every cHSS $\mathcal{S}=G/P$ is associated to a marked Dynkin diagram $(\mathcal{D}(G),\gamma)$ with $\gamma$ being a long root. The marking at $\gamma$ gives rise to a partition of set of roots $\Delta=\Delta^-\sqcup\Delta^0\sqcup\Delta^+$, corresponding to compact roots, negative and positive non-compact roots respectively. Denote by $E^+$ the positive non-compact root vectors, then span$\{E^+\}=\frak{m}^+\cong T_o\mathcal{S}$. Furthermore, we have a finer disjoint partition $E^+=\{E_{\gamma}\}\sqcup H_{\gamma}\sqcup N_{\gamma}$ with $\alpha=E_{\gamma}, \mathscr{H}_{\alpha}=$span$\{H_{\gamma}\}, \mathscr{N}_{\alpha}=$span$\{N_{\gamma}\}$. Let $\theta$ denote the unique node which is adjacent to the marked root $\gamma$ in the Dynkin diagram if $(\mathcal{D},\gamma)$ is not associated to Grassmannian of rank $\geq 2$ and let $\theta_1,\theta_2$ denote the two nodes adjacent to $\gamma$ if $(\mathcal{D}(G),\gamma)$ is associated to Grassmannian of rank $\geq 2$. Then the following lemma follows immediately from Proposition \ref{Prop 1.8 for paper 2},

\begin{lemma}
\label{Lemma 2 for paper 2}
If $(\mathcal{D}(G),\gamma)$ is not associated to Grassmannian of rank $\geq 2$, then $H_{\gamma}=\{E_{\gamma+\theta+\cdots}: \textnormal{if}\ \ \gamma+\theta+\cdots$ is a root\}, $N_{\gamma}=\{E_{\gamma+2\theta+\cdots}:$ if $\gamma+2\theta+\cdots$ is a root\}. Otherwise, $H_{\gamma}=\{E_{\gamma+\theta_1+\cdots}:$ if $\gamma+\theta_1+\cdots$ is a root\}$\sqcup \{E_{\gamma+\theta_2+\cdots}:$ if $\gamma+\theta_2+\cdots$ is a root\}, $N_{\gamma}=\{E_{\gamma+\theta_1+\theta_2+\cdots}:$ if $\gamma+\theta_1+\theta_2+\cdots$ is a root\}.
\end{lemma} 

Denote by $(\mathbb{C}^m\cong\frak{m}^+,z_i)$ the coordinate dual to the positive root vectors in $E^+$, by virtue of Harish-Chandra embedding (cf. chapter 5, [Mok 89]), it serves as a Harish-Chandra coordinate chart. The following table lists all the cHSS and its associated marked Dynkin diagram $(\mathcal{D}(G),\gamma)$ excluding $Q^1,Q^2$ which are isomorphic to degree 2 rational curves in $\mathbb{P}^2$ and reducible $\mathbb{P}^1\times\mathbb{P}^1$ respectively.
\newpage
\begin{table}
\caption{Hermitian Symmetric Spaces and marked Dynkin diagrams cf. \cite{[Hel 78]}}
\label{tab:2}       
\begin{tabular}{|l|l|}
	\hline
	Type & Marked Dynkin diagram\\
	\hline
	$G(p,q),  p,q\geq 1$ & $A_{p+q-1}:$ 
	\xymatrix@C=1em@R=1em{
\circ_1\ar@{-}[r]&\cdots\ar@{-}[r]&\circ_{p-1}\ar@{-}[r]&\bullet_{p}\ar@{-}[r]&\circ_{p+1}\ar@{-}[r]&\cdots\ar@{-}[r]&\circ_{p+q-1}
} \\
  \hline
	$G^{\Rmnum{2}}(n,n), n\geq 3$ & $D_n:$
	\xymatrix@C=1em@R=1em{
&\circ_{n-1}\ar@{-}[d] \\
\bullet_n \ar@{-}[r]&\circ_{n-2}\ar@{-}[r]&\circ_{n-3}\ar@{-}[r]&\cdots\ar@{-}[r]&\circ_2\ar@{-}[r]&\circ_1
} \\
	\hline
	$G^{\Rmnum{3}}(n,n), n\geq 2$ & $C_n:$
	\xymatrix@C=1em@R=1em{
\circ_1\ar@{-}[r]&\circ_2\ar@{-}[r]&\cdots\ar@{-}[r]&\circ_{n-1}\ar@{<=}[r]&\bullet_n
}\\
	\hline
	$Q^{2m-2}, m\geq 3$ & $D_m:$
	\xymatrix@C=1em@R=1em{
&\circ_{m-1}\ar@{-}[d] \\
\circ_m \ar@{-}[r]&\circ_{m-2}\ar@{-}[r]&\circ_{m-3}\ar@{-}[r]&\cdots\ar@{-}[r]&\circ_2\ar@{-}[r]&\bullet_1
} \\
	\hline
	$Q^{2m-1}, m\geq 2$ & $B_m:$
	\xymatrix@C=1em@R=1em{
\bullet_1\ar@{-}[r]&\circ_2\ar@{-}[r]&\cdots\ar@{-}[r]&\circ_{m-1}\ar@{=>}[r]&\circ_m
} \\
	
\hline

	$\Rmnum{5}$ & $E_6:$
	\xymatrix@C=1em@R=1em{
&&\circ_{2}\ar@{-}[d] \\
\bullet_{6} \ar@{-}[r]&\circ_{5}\ar@{-}[r]&\circ_{4}\ar@{-}[r]&\circ_{3}\ar@{-}[r]&\circ_{1}
}\\

\hline

$\Rmnum{6}$ & $E_7$
	\xymatrix@C=1em@R=1em{
&&&\circ_{2}\ar@{-}[d] \\
\bullet_{7} \ar@{-}[r]&\circ_{6}\ar@{-}[r]&\circ_{5}\ar@{-}[r]&\circ_{4}\ar@{-}[r]&\circ_{3}\ar@{-}[r]&\circ_{1}
} \\
	\hline
	\end{tabular}
	
	\end{table}

We try to construct equivariant embedding $i:\mathcal{S}_0\hookrightarrow\mathcal{S}$ from isomorphism between their root systems, the reader may refer to \cite{[Hum 72]} for relevant theories concerning Lie algebra and its root system. In particular, we introduce the notion of \textsl{root correspondence} as follows:

\begin{definition}
\label{def 6 for paper 2}
Let $\mathcal{S}_0,\mathcal{S}$ be cHSS associated to marked Dynkin diagrams $(\mathcal{D}(G_0),\gamma_0)$, $(\mathcal{D}(G),\gamma)$ respectively and $\Delta_0,\Delta$ denote the respective root systems for $G_0$=Aut$(\mathcal{S}_0)$, $G=$Aut$(\mathcal{S})$ respectively. An injective mapping $\Phi:\Delta_0\rightarrow\Delta$ is said to be a root correspondence if $\Phi$ is an isomorphism onto its image and $\Phi(\gamma_0)=\gamma$, i.e., preserving marked root.
\end{definition}Note since $\Phi$ is an isomorphism from root system onto its image, it induces an injective homomorphism between Lie algebras $\widetilde{\Phi}:\frak{g}_0\hookrightarrow\frak{g}$, i.e., $\widetilde{\Phi}$ is an isomorphism from Lie algebra $\frak{g}_0$ onto its image. 

Given an arbitrary pair $(\mathcal{S}_0,\mathcal{S})$, such root correspondence may not always exist. But if the pair is of sub-diagram type (cf. Theorem \ref{Thm 1.1} ), then $\Phi$ is just obtained by identification of the nodes when $\mathcal{D}(G_0)$ is realized as a subdiagram of $\mathcal{D}(G)$ with marked root $\gamma_0$ being simultaneously identified with marked root $\gamma$.

\begin{lemma}
\label{Lemma 3 for paper 2}
	Let $(\mathcal{S}_0,\mathcal{S})$ be a pair of cHSS associated to marked Dynkin diagrams $(\mathcal{D}(G_0),\gamma_0)$,$(\mathcal{D}(G),\gamma)$ respectively, and suppose there exists a root correspondence $\Phi:\Delta_0\rightarrow\Delta$, then  
	\begin{enumerate}
	\item $\Phi(\Delta_0^+)\subset\Delta^+,\ \ \Phi(\Delta_0^-)\subset\Delta^-,\ \ \Phi(\Delta_0^0)\subset\Delta^0$;
	\item Let $\widetilde{\Phi}:\frak{g}_0\hookrightarrow\frak{g}$ be the isomorphism from $\frak{g}_0$ onto its image induced by $\Phi$, recall we have the partition $E_0^+=\{\alpha_0=E_{\gamma_0}\}\sqcup H_{\gamma_0}\sqcup N_{\gamma_0},\ \ E^+=\{\alpha=E_{\gamma}\}\sqcup H_{\gamma}\sqcup N_{\gamma}$ for $\mathcal{S}_0,\mathcal{S}$ respectively, then $\widetilde{\Phi}(\alpha_0)=\alpha,\ \ \widetilde{\Phi}(H_{\gamma_0})\subset H_{\gamma},\ \ \widetilde{\Phi}(N_{\gamma_0})\subset N_{\gamma}$.
	\end{enumerate}
\end{lemma}

\textbf{Proof} Since $\phi$ is an isomorphism onto its image of root systems, let $\alpha,\beta\in\Delta_0$, then $\alpha+\beta\in\Delta_0$ if and only if$\phi(\alpha+\beta)\in\Delta$. In particular, for simple roots $\gamma_1,\gamma_2,\ \ \gamma_1-\gamma_2$ is never a root. Then this lemma quickly follows from this property and the description of root vectors in $H_{\gamma},N_{\gamma}$ as shown in Lemma \ref{Lemma 2 for paper 2}

\begin{flushright}
QED
\end{flushright}

The existence of such root correspondence has significance in judging whether $(\mathcal{S}_0,\mathcal{S})$ is an admissible pair.

\begin{proposition}
\label{Prop 2 for paper 2}
The root correspondence $\Phi:\Delta_0\rightarrow\Delta$ naturally induces a standard embedding $i_{\Phi}:\mathcal{S}_0\hookrightarrow\mathcal{S}$, making $(\mathcal{S}_0,\mathcal{S})$ an admissible pair.
\end{proposition}

\textbf{Proof} $\Phi$ gives rise to an injective homomorphism $\widetilde{\Phi}:\frak{g}_0\hookrightarrow\frak{g}$. By Lemma \ref{Lemma 3 for paper 2}, restriction of $\widetilde{\Phi}$ to $\frak{m}^-_0,\frak{k^{\mathbb{C}}}_0,\frak{m}^+_0$ are all injective homomorphisms on the three Lie subalgebras arising from Harish-Chandra decomposition, i.e., $\frak{g}_0=\frak{m}^-_0+\frak{k^{\mathbb{C}}}_0+\frak{m}^+_0,\ \ \frak{g}=\frak{m}^-+\frak{k^{\mathbb{C}}}+\frak{m}^+$ and $\widetilde{\Phi}(\frak{m}^{-}_0)\subset\frak{m}^{-},\widetilde{\Phi}(\frak{k}^{\mathbb{C}}_0)\subset\frak{k}^{\mathbb{C}},\widetilde{\Phi}(\frak{m}^+_0)\subset\frak{m}^+$. Then $\widetilde{\Phi}$ induces an equivariant embedding $i_{\Phi}:\mathcal{S}_0\hookrightarrow\mathcal{S}$. In particular, on some Harish-Chandra coordinate chart $\mathbb{C}^n\cong U_0\subset\mathcal{S}_0,\mathbb{C}^m\cong U\subset\mathcal{S}$, the linear mapping $\widetilde{\Phi}:\frak{m}^+_0\hookrightarrow\frak{m}^+$ induces $i_{\Phi}(U_0)\subset U$. So $U_0$ is realized as an affine subspace of $U$, in particularly it maps affine lines in $U_0$ to affine lines in $U$, which means nothing but mapping minimal rational curves of $\mathcal{S}_0$ to those of $\mathcal{S}$. If $\{z_i\}^n_{i=1},\ \ \{\omega_j\}^m_{j=1}$ denotes the coordinate on $U_0,U$ dual to positive root vectors $E^+_0,E^+$ respectively, then we may express the embedding as

\begin{align*}
w_i & =z_i,\ \ i=1,2,...,n\\
w_j & =0,\ \ j=n+1,...,m
\end{align*} To complete the proof, we only need to show $i_{\Phi}$ respects VMRT.
	
	Fix any reference point $o$, by abuse of notation, we do not differentiate between $o$ and its image $i(o)$. We have identifications $T_o(\mathcal{S}_0)\cong\frak{m}^+_0,\ \ T_{o}\mathcal{S}\cong\frak{m}^+$, as a result $\forall v\in T_o(\mathcal{S}_0)$, the differential $di_{\Phi}(v)=\widetilde{\Phi}(v)$. For any $\beta\in\widetilde{\mathscr{C}}_o(\mathcal{S}_0)$, there exists some positive root vector $\alpha\in E^+_0$ and $g\in\frak{k}_{s0}=[\frak{k}^{\mathbb{C}}_0,\frak{k}^{\mathbb{C}}_0]$, such that $\beta=$ exp(ad$g)(\alpha)$, so $di_{\Phi}(\beta)= \widetilde{\Phi}($exp(ad$g)(\alpha)$) = exp(ad$\widetilde{\Phi}(g))(\widetilde{\Phi}(\alpha))$. Note that $\widetilde{\Phi}(\alpha)\subset E^+$ is again some positive root vector, and $\widetilde{\Phi}(g)\subset\frak{k}_s$, so $di_{\Phi}(\beta)\subset\widetilde{\mathscr{C}}_o(\mathcal{S})$. On the other hand, by Lemma \ref{Lemma 3 for paper 2}, $di_{\Phi}(\mathscr{H}_{\alpha_0})\subset\mathscr{H}_{\alpha},\ \ di_{\Phi}(\mathscr{N}_{\alpha_0})\subset \mathscr{N}_{\alpha}$, and the fact that the ambient VMRT $\mathscr{C}_o(\mathcal{S})$ is the Zariski closure of some quadratic polynomial $\zeta(\xi)$ for $\xi\in\mathscr{H}_{\alpha}$, then the sub-VMRT $\mathscr{C}_o(\mathcal{S}_0)\subset\mathscr{C}_o(\mathcal{S})$ is obtained by restriction of $\zeta$ to the subspace $\mathscr{H}_{\alpha_0}\subset\mathscr{H}_{\alpha}$, so $di_{\Phi}$ must respects VMRT.

\begin{flushright}
QED
\end{flushright}


\section {Classification of admissible pairs of irreducible compact Hermitian Symmetric Spaces}
\label{chapter 4}


\subsection{\bf{Candidates of admissible pairs of cHSS}}
\label{section 4.1}

In this section, we first of all work out a list of all the possible admssible pairs $(\mathcal{S}_0,\mathcal{S})$ of cHSS as our Main Theorem 1 shows. We assume $\mathcal{S}_0$ to be non-linear throughout this section. The following lemma is a necessary condition for $(\mathcal{S}_0,\mathcal{S})$ to be an admissible pair which is immediate from the definition.

\begin{lemma}
\label{Lemma 4.1}
If $(\mathcal{S}_0,\mathcal{S})$\ \ is admissible with dim$(\mathcal{S}_0)=n<$dim$(\mathcal{S})=m$, then for their respective (embedded) VMRTs $\mathscr{C}_o(\mathcal{S}_0)\subset\mathbb{P}T_o(\mathcal{S}_0)\cong\mathbb{P}^{n-1}$,$\mathscr{C}_o(\mathcal{S})\subset\mathbb{P}T_o(\mathcal{S})\cong\mathbb{P}^{m-1}$, there exists some projective subspace $\mathbb{P}(V)\cong\mathbb{P}T_o(\mathcal{S}_0)$ such that $\mathscr{C}_o(\mathcal{S}_0)\subset\mathbb{P}T_o(\mathcal{S}_0)$ is projectively equivalent to $\mathscr{C}_o(\mathcal{S})\cap\mathbb{P}(V)\subset\mathbb{P}(V)$.
\end{lemma}

In other words, there exists an embedding between abstract VMRTs $\psi:\mathcal{A}_0\hookrightarrow\mathcal{A}$ such that $\phi_0=\phi\circ\psi$ where $\phi_0:\mathcal{A}_0\hookrightarrow\mathbb{P}^{n-1}, \phi:\mathcal{A}\hookrightarrow\mathbb{P}^{m-1}$ are the embeddings whose images are the embedded VMRTs $\mathscr{C}_o(\mathcal{S}_0), \mathscr{C}_o(\mathcal{S})$ respectively. In the sense of the above Lemma, we also say that $(\mathcal{A}_0,\mathcal{A})$ is an \textsl{admissible pair of VMRT}.


\begin{lemma}\label{Lemma 4.2}
Suppose $(\mathcal{S}_0,\mathcal{S})$ is an admissible pair of cHSS,

\begin{enumerate}
\item if $\mathcal{S}$ is a type $\Rmnum{1}$ Grassmannian $G(r,s)$, then $\mathcal{S}_0$ is either a type $\Rmnum{1}$ Grassmannian $G(p,q)$ as well with $2\leq p\leq r,\ \ 2\leq p\leq s$ (this pair is of sub-diagram type, in category 1), or it is a type $\Rmnum{3}$ Grassmannian $G^{\Rmnum{3}}(n,n)$ with $3\leq n\leq$min\{$r,s$\}. Note $G^{\Rmnum{3}}(2,2)\cong Q^3$ (this pair is of special type in category 3);

\item if $\mathcal{S}$ is hyperquadric $Q^m$, then $\mathcal{S}_0$ can only be hyperquadric $Q^n$ as well with $n<m$ ( this pair is of deletion type in category 2 or of special type in category 3 depending on $m-n$ is even or odd number respectively);

\item if $\mathcal{S}$ is Lagrangian Grassmannian, then so is $\mathcal{S}_0$. On the other hand, if $\mathcal{S}_0$ is symplectic Grassmannian, then Lagrangian Grassmannian is the only possibility for $\mathcal{S}$ except the usual Grassmannian as described in 1 (this pair is of sub-diagram type in category 1).

\end{enumerate}
\end{lemma}

\textbf{proof} For 1, the abstract VMRT $\mathcal{A}$ of $G(r,s)$ is $\mathbb{P}^{r-1}\times\mathbb{P}^{s-1}$ and Segre embedding $\phi:\mathbb{P}^{r-1}\times\mathbb{P}^{s-1}\hookrightarrow\mathbb{P}^{rs-1}$ realizes $\mathcal{A}$ as embedded VMRT of $G(r,s)$. If $\mathcal{A}_0$ is not abstract VMRT of any projective space, i.e., $\mathcal{A}_0$ is non linear, then for $\mathcal{A}_0\stackrel{\psi}{\hookrightarrow}\mathcal{A}\stackrel{\phi}{\hookrightarrow}\mathbb{P}^{rs-1}$, either $\mathcal{A}_0$ is reducible and $\phi_0=\phi\circ\psi$ is realized by $\mathcal{O}(1)$ or $\mathcal{A}_0$ is irreducible and $\phi_0$ is realized by $\mathcal{O}(2)$, which  correspond to $\mathcal{A}_0=\mathbb{P}^{p-1}\times\mathbb{P}^{q-1},\ \ \mathcal{S}_0=G(p,q),\ \ (2\leq p<r,2\leq q<s)$ and $\mathcal{A}_0=\mathbb{P}^{n-1},\ \ \mathcal{S}_0=G^{\Rmnum{3}}(n,n),\ \ 2\leq n\leq$min$\{r,s\}$.

For 2, $\mathcal{A}=Q^{m-2}\cap\mathbb{P}^{n-1}$ is a hypersurface in $\mathbb{P}^{n-1},\ \ n=$dim$(\mathcal{S}_0)$, the reduced cycle $\mathcal{A}_0\subset\mathbb{P}^{n-1}$ is of degree 1 or 2, hence $\mathbb{P}^k$ or $Q^l$. Since the case with $\mathcal{A}_0$ being linear is excluded, then $\mathcal{A}_0$ can only be another hyperquadric, hence $\mathcal{S}_0$ is hyperquadric. 

For 3, it is because that the VMRT of $G^{\Rmnum{3}}(m,m)$ is $\mathbb{P}^{m-1}\hookrightarrow\mathbb{P}^{m(m+1)/2-1}$ embedded through linear system of $\mathcal{O}(2)$ while all the other embedded VMRT are embedded through $\mathcal{O}(1)$. 

\begin{flushright}
QED
\end{flushright}

Lemma \ref{Lemma 4.2} exhausts all the admissible pairs $(\mathcal{S}_0,\mathcal{S})$ in which Lagrangian Grassmannian $G^{\Rmnum{3}}(m,m)$ is involved. In the rest of this section, we are justified to exclude Lagrangian Grassmannian whose embedded VMRT contains no projective line.


\begin{lemma}
\label{Lemma 4.3}
Let $(\mathcal{A}_0,\mathcal{A})$ be an admissible pair of abstract VMRT, for their respective embedded VMRTs $\mathscr{C}_o(\mathcal{S}_0),\mathscr{C}_o(\mathcal{S})$, there exists some $\mathbb{P}(T_o\mathcal{S}_0)\cong\mathbb{P}(V)\subset\mathbb{P}(T_o\mathcal{S})$, such that $\mathscr{C}_o(\mathcal{S}_0)=\mathbb{P}(V)\cap\mathscr{C}_o(\mathcal{S})\hookrightarrow\mathscr{C}_o(\mathcal{S})$ is an isometric totally geodesic embedding.    
\end{lemma}

\textbf{Proof} Since the inclusion embedding is equivariant, it is naturally isometric with respect to some fixed canonical metric $h_0,h$ on $\mathcal{A}_0,\mathcal{A}$ respectively. 

For simplicity, we write $\mathbb{P}^m$ for $\mathbb{P}(T_o\mathcal{S}),\ \ \mathbb{P}^n$ for $\mathbb{P}(V)$. We choose some $\mathbb{P}^2\subset\mathbb{P}^m$, but $\not\subset\mathscr{C}_o(\mathcal{S})$ such that $C=\mathbb{P}^2\cap\mathscr{C}_o(\mathcal{S})$ is a degree 2 rational curve. Note that every such curve $C\subset\mathscr{C}_o(\mathcal{S})$ is obtained in this way. Since the normal bundle $N_{C|\mathscr{C}_o(\mathcal{S})}=T_{\mathscr{C}_o(\mathcal{S})}|_C/T_C$ is a proper sub-bundle of $T_{\mathbb{P}^m}|_C/T_{\mathbb{P}^2}|_C\cong\mathcal{O}(2)^{m-2}$, then positivity of $T_{\mathscr{C}_o(\mathcal{S})}$ yields the splitting type $N_{C|\mathscr{C}_o(\mathcal{S})}=\mathcal{O}(2)^a\oplus\mathcal{O}(1)^b\oplus\mathcal{O}^c$. Now we show this splitting type forces the vanishing of second fundamental form $\sigma_{C|\mathscr{C}_o(\mathcal{S})}$.



In fact, $\sigma_{C|\mathscr{C}_o(\mathcal{S})}$ is topologically isomorphic to $N^*_{C|\mathscr{C}_o(\mathcal{S})}\otimes T_C$-valued $(0,1)$-form on $C$ (cf. \cite{[Mok 05]}), one of whose global sections is denoted by $\eta$, and we decompose $\eta$ into line bundle valued $(0,1)$-form according to the splitting type $N^*_{C|\mathscr{C}_o(\mathcal{S})}\otimes T_C=\mathcal{O}(2)^c\oplus\mathcal{O}(1)^b\oplus\mathcal{O}^a$, i.e., $\eta=\eta_1+\cdots+\eta_k$. If $\sigma_{C|\mathscr{C}_o(\mathcal{S})}$ does not vanish at some $p\in C$, then so is at least one of the components, say $\eta_1$. Since all the vector bundles involved here is homogeneous with respect to isometry group, then we produce a nowhere vanishing section of line bundle $\mathcal{O}(l)$-valued $(0,1)$-form ($l=0,1,2$). By means of canonical metric, $(0,1)$-form is topologically isomorphic to $T_C$. As a result, the nowhere vanishing $\eta_1$ is a smooth section of line bundle $\mathcal{O}(l+2)$, whose topological class never allows the existence of nowhere vanishing section, thus we get to contradiction. So $\sigma_{C|\mathscr{C}_o(\mathcal{S})}$ must vanish everywhere, i.e., $C\subset\mathscr{C}_o(\mathcal{S})$ is totally geodesic.

Having proven the total geodesy of any degree 2 rational curve in $\mathscr{C}_o(\mathcal{S})$, we consider this lemma in its full strength. Let $\phi_0:\mathcal{A}_0\cong\mathscr{C}_o(\mathcal{S}_0)=\mathbb{P}^n\cap\mathscr{C}_o(\mathcal{S})\subset\mathbb{P}^n\subset\mathbb{P}^m$ be the embedding realizing $\mathcal{A}_0$ as embedded VMRT. Now fix any point $p\in\mathscr{C}_o(\mathcal{S}_0)$ and consider a degree 2 rational curve $C\subset\mathscr{C}_o(\mathcal{S}_0)\subset \mathscr{C}_o(\mathcal{S})$ passing through $p$, we furthermore assume $C$ is obtained through:

\begin{center}
$\mathbb{P}^1\stackrel{d}{\hookrightarrow}\mathbb{P}^1\times\mathbb{P}^1\stackrel{i}{\hookrightarrow}\mathcal{A}_0\stackrel{\phi_0}{\cong}\mathscr{C}_o(\mathcal{S}_0)$
\end{center}where $d$ is the diagonal embedding and $i$ realizes the product $\mathbb{P}^1\times\mathbb{P}^1$ as a totally geodesic submanifold of $\mathcal{A}_0$ by virtue of polysphere theorem (cf. chapter 5, [Mok 89]) such that each factor $\phi_0\circ i(\mathbb{P}^1\times \{o\}),\ \ \phi_0\circ i(\{o\}\times \mathbb{P}^1)$ is a minimal rational curve on $\mathscr{C}_o(\mathcal{S})$ (projective line in $\mathbb{P}^m$) passing through $p$. Denote by $v_1,v_2$ the tangent vector at $p$ of the first and second factor respectively. Then by total geodesy of the degree 2 rational curve $\phi_0\circ i\circ d(\mathbb{P}^1)=C\subset\mathscr{C}_o(\mathcal{S})$, we have 

\begin{center}
$\sigma(v_1+v_2,v_1+v_2)=0\Rightarrow \sigma(v_1,v_2)=0$ 
\end{center}since $\sigma(v_1,v_1)=\sigma(v_2,v_2)=0$ for they are minimal rational tangents, where $\sigma$ is the second fundamental form of $\mathscr{C}_o(\mathcal{S}_0)$ as a submanifold of $\mathscr{C}_o(\mathcal{S})$. On the other hand, it is well known from the structure theory of VMRT that $v_2\in\mathscr{N}_{v_1}\cap T_p(\mathscr{C}_o(\mathcal{S}_0))$ where $\mathscr{N}_{v_1}$ is in the sense of the decomposition $T_p(\mathscr{C}_o(\mathcal{S}))=\mathbb{C}v_1\oplus\mathscr{H}_{v_1}\oplus\mathscr{N}_{v_1}$ (cf. chapter 7, [Mok 89]). Having achieved this, we can just use the polarization arguments to obtain $\sigma(v,u)=0$ for any $v,u\in T_p(\mathscr{C}_o(\mathcal{S}_0))$. The procedure is, roughly speaking, varying $v_1(t)$ holomorphically with respect to $t$ while the corresponding $v_2(t)\in\mathscr{N}_{v_1(t)}$ will change anti-holomorphically w.r.t. $t$, the reader may refer to proposition 3 on page 111 in [Mok 89] for more details. Hence we have established the total geodesy of $\mathscr{C}_o(\mathcal{S}_0)\subset\mathscr{C}_o(\mathcal{S})$.

\begin{flushright}
QED
\end{flushright}

\begin{proposition}
\label{Prop 4.3}
Suppose $\mathcal{S}$ (dim$(\mathcal{S})=m$) is a cHSS with its embedded VMRT $\mathscr{C}_o(\mathcal{S})$ at some fixed reference point $o\in\mathcal{S}$. Let $V\subset T_o\mathcal{S}$ be an $n$-dimensional vector subspace such that $\mathscr{C}_o(\mathcal{S})\cap\mathbb{P}(V)$ is projectively equivalent to the embedded VMRT $\mathscr{C}_o(\mathcal{S}_0)\subset\mathbb{P}^{n-1}$ of another cHSS $\mathcal{S}_0$ of dimension $n$, then $(\mathcal{S}_0,\mathcal{S})$ is an admissible pair, i.e., $\mathcal{S}_0\subset\mathcal{S}$ as a totally geodesic standard model can be uniquely constructed out of the datum $V,\mathscr{C}_o(\mathcal{S}_0),\mathcal{S}$.
\end{proposition}


\textbf{Proof}  Let $\frak{g}=\frak{m}^-\oplus\frak{k}^{\mathbb{C}}\oplus\frak{m}^+$ be the Harish-Chandra decomposition which is uniquely associated to cHSS $\mathcal{S}$ (cf. [Mok 89]). We only need to recover three subalgebras $\frak{m}_0^-\subset\frak{m}^-,\frak{k}_0^{\mathbb{C}}\subset\frak{k}^{\mathbb{C}},\frak{m}_0^+\subset\frak{m}^+$ such that $\frak{g}^{\mathbb{C}}_0=\frak{m}_0^-\oplus\frak{k}_0^{\mathbb{C}}\oplus\frak{m}_0^+$ is the Harish-Chandra decomposition of some Lie subalgebra $\frak{g}^{\mathbb{C}}_0\subset\frak{g}^{\mathbb{C}}$. 

$\frak{m}_0^+$ is just identified with the linear subspace $V$, making it an Abelian subalgebra of $\frak{m}^+\cong\mathbb{C}^m$. Through the embedding $\mathscr{C}_o(\mathcal{S}_0)=\mathscr{C}_o(\mathcal{S})\cap\mathbb{P}(V)\hookrightarrow\mathscr{C}_o(\mathcal{S})$, we can recover a complex Lie subalgebra $\frak{k}^{\mathbb{C}}_0\subset\frak{k}^{\mathbb{C}}$ induced by their respective real forms $\frak{k}_0\subset\frak{k}$ such that the center $z\in\frak{k}$ is also contained in $\frak{k}_0$.

Now we define $\frak{m}_0^-:=\overline{V}\subset\frak{m}^-$, where bar $\overline{\cdot}$ is the complex conjugation, we have the following claim:

\begin{center}
$[[\overline{\eta},\alpha],\beta]\in V\cong\frak{m}^+_0$\ \ for $\forall\eta,\alpha,\beta\in \frak{m}_0^+\subset\frak{m}^+$\ \ \ \ \ \ (\dag)
\end{center}

Assuming the claim $(\dag)$,\ \ we show that $\frak{m}_0^-:=\overline{V}$\ \ is what we desire, i.e., $\frak{g}^{\mathbb{C}}_0:=\frak{m}_0^+\oplus\frak{k}_0^{\mathbb{C}}\oplus\frak{m}_0^-$\ \ is a Lie subalgebra of $\frak{g}^{\mathbb{C}}$\ \ with the Harish-Chandra decomposition.

For $\forall\eta, \beta\in\frak{m}_0^+$,\ \ we show that $[\overline{\eta},\beta]\in\frak{k}_0^{\mathbb{C}}$.\ \ Take any $\alpha\in\widetilde{\mathcal{A}_0}$,\ \ since $[\overline{\eta},\beta]\in\frak{k}^{\mathbb{C}}=\frak{aut}(\widetilde{\mathcal{A}})$,\ \ we have $[[\overline{\eta},\beta],\alpha]\in T_{[\alpha]}(\widetilde{\mathcal{A}})$\ \ (cf. Lemma 1.1.1 [HwM 05]), then by the claim $(\dag),\ \ [[\overline{\eta},\beta],\alpha]\in T_{[\alpha]}(\widetilde{\mathcal{A}_0})$,\ \ so $[\overline{\eta},\beta]\in\frak{aut}(\widetilde{\mathcal{A}_0})=\frak{k}_0^{\mathbb{C}}$.\ \ We showed $[\frak{m}_0^+,\frak{m}_0^-]\subset\frak{k}_0^{\mathbb{C}}$.

For $\xi\in\frak{k}_0^{\mathbb{C}}$, since $\frak{k}_0^{\mathbb{C}}=\frak{aut}(\widetilde{\mathcal{A}_0})$\ \ and $\mathcal{A}_0$\ \ is a non-degenerate subvariety, we get $[\frak{k}_0^{\mathbb{C}},\frak{m}_0^+]\subset\frak{m}_0^+$.\ \ It follows that $[\xi,\overline{\eta}]]=\overline{[\overline{\xi},\eta]}\in\overline{\frak{m}_0^+}=\frak{m}_0^-$.

Now we come to the proof of the claim $(\dag)$. Fixing some canonical metric $h$ on $\mathcal{S}$ and $\alpha\in\widetilde{\mathscr{C}}_o(\mathcal{S})\cap V=\widetilde{\mathscr{C}}_o(\mathcal{S}_0)$. As a result of total geodesy thanks to Lemma \ref{Lemma 4.3}, the two second fundamental forms of $\widetilde{\mathscr{C}}_o(\mathcal{S}_0),\widetilde{\mathscr{C}}_o(\mathcal{S})$ as submanifolds in the ambient space $T_o\mathcal{S}$ w.r.t. Euclidean flat connection, denoted by $\sigma'_0,\sigma'$ respectively, are identical in their common domain of definition $\mathscr{H}_{\alpha,0}=\mathscr{H}_{\alpha}\cap V$, where $\mathscr{H}_{\alpha}=T_{[\alpha]}(\mathscr{C}_o(\mathcal{S}))$ is in the sense of eigenspace decomposition $T_o\mathcal{S}=\mathbb{C}\alpha\oplus\mathscr{H}_{\alpha}\oplus\mathscr{N}_{\alpha}$ hence $\mathscr{H}_{\alpha,0}=T_{[\alpha]}(\mathscr{C}_o(\mathcal{S}_0))$. Denote by $\mathscr{N}_{\alpha,0}$ the normal complement subspace of $\mathbb{C}\alpha\oplus\mathscr{H}_{\alpha,0}$ in $V$, then $\mathscr{N}_{\alpha,0}$ is generated by taking second fundamental form of vectors in $\mathscr{H}_{\alpha,0}$, i.e., $\mathscr{N}_{\alpha,0}=\{\sigma'_0(\xi_1,\xi_2):\xi_1,\xi_2\in\mathscr{H}_{\alpha,0}\}$, for $\mathscr{C}_o(\mathcal{S}_0)$ is non-degenerate in $\mathbb{P}(V)$, i.e., $\mathscr{C}_o(\mathcal{S}_0)\subset\mathbb{P}(V)$ is not contained in any hyperplane of $\mathbb{P}(V)$ and $\mathscr{C}_o(\mathcal{S}_0)$ can be expressed as the Zariski closure of a vector valued polynomial $\zeta(\xi)=\alpha+\xi+\sigma'(\xi,\xi)$ (cf. Lemma \ref{Lemma 1 of paper 2}).  Likewise, $\mathscr{N}_{\alpha}=\{\sigma'(\xi_1,\xi_2):\xi_1,\xi_2\in\mathscr{H}_{\alpha}\}$. As a consequence, we are led to an ideal decomposition on $V$ as follows:

\begin{center}
$V=\mathbb{C}\alpha\oplus\mathscr{H}_{\alpha,0}\oplus\mathscr{N}_{\alpha,0}$,
\end{center}such that $\mathscr{N}_{\alpha,0}=\mathscr{N}_{\alpha}\cap V$. Let $X,Y,Z\in\frak{m}^+\cong T_o(\mathcal{S})$, the ralation between Riemannian curvature tensor $R$ and the Lie bracket is $R(X,Y)Z=[[X,Y],Z]$ (cf. [Mok 89]). Then we again borrow the polorization argument from [Mok 89] to obtain the desired result $R(\xi,\overline{\alpha})(\beta)=[[\xi,\overline{\alpha}],\beta]\in V,\ \ \forall\alpha,\beta,\xi\in V$, thus we complete the proof.

\begin{flushright}
QED
\end{flushright}

This proposition allows us to do ``reduction'' on admissible pairs, i.e., if $(\mathcal{A}_0,\mathcal{A})$ is an admissible pair of VMRT, then $(\mathcal{S}_0,\mathcal{S})$ is an admissible pair in the sense of Def \ref{Def 1.2 for paper 2} and vice versa.

\begin{lemma}
\label{Lemma 4.4}
Suppose $(\mathcal{S}_0,\mathcal{S})$\ \ is an admissible pair and $\mathcal{S}_0$\ \ is an exceptional Hermitian symmetric space, i.e., either $\Rmnum{5}=E_6/(\textnormal{Spin}(10)\times U(1))$\ \ or $\Rmnum{6}=E_7/(E_6\times U(1))$,\ \ then $\mathcal{S}$\ \ must also be an exceptional Hermitian symmetric space, i.e., the only admissible pair for $\mathcal{S}_0$\ \ to be an exceptional space is $(\Rmnum{5},\Rmnum{6})$.
\end{lemma}

\textbf{Proof} Since VMRT of $\Rmnum{6}$ is $\Rmnum{5}$, and VMRT of the latter is $G^{\Rmnum{2}}(5,5)$, we only need to determine admissible pairs $(G^{\Rmnum{2}}(5,5),\mathcal{A})$ where $\mathcal{A}$ is itself a VMRT of another cHSS. By virtue of Lemma \ref{Lemma 4.2}, we can exclude the cases in which $\mathcal{A}$ is VMRT of Grassmannian, hyperquadric and symplectic Grassmannian. Since $G^{\Rmnum{2}}(n,n)$ is not a VMRT for $n\geq 6$, nor is $\Rmnum{6}$, the only remaining possibility is $\mathcal{A}=\Rmnum{5}$, and $(G^{\Rmnum{2}}(5,5),\Rmnum{5})$ is actually admissible which means, thanks to Prop \ref{Prop 4.3}, $(\Rmnum{5},\Rmnum{6})$ is also an admissible pair (this pair is of deletion type in category 2).   

\begin{flushright}
QED
\end{flushright}

\begin{lemma}
\label{Lemma 4.5}
Suppose $(\mathcal{S}_0,G^{\Rmnum{2}}(m,m))\ \ m\geq 5$\ \ is an admissible pair, then $\mathcal{S}_0$\ \ falls into one of the following three cases:

\begin{enumerate}
\item $\mathcal{S}_0=G(r,s)$ where $2\leq r ,s$ and $r+s\leq m$ (if either $r$ or $s$ is 2, then $(G(r,s),G^{\Rmnum{2}}(m,m))$ is of deletion type in category 2; if both $r,s\geq 3$, then the pair is of special type in category 3).

\item $\mathcal{S}_0=Q^n,\ \ n=2,3,...,6$. ($(Q^6,G^{\Rmnum{2}}(m,m))$ is of sub-diagram type in category 1, note that $G^{\Rmnum{2}}(4,4)\cong Q^6$ and all the rest are in category 4).

\item $\mathcal{S}_0=G^{\Rmnum{2}}(n,n)\ \ 4\leq n<m$ (this is of sub-diagram type, note that $G^{\Rmnum{2}}(2,2)\cong\mathbb{P}^1,\ \ G^{\Rmnum{2}}(3,3)\cong\mathbb{P}^3$).
\end{enumerate}

\end{lemma}

\textbf{Proof} According to 3 of Lemma \ref{Lemma 4.2} and Lemma \ref{Lemma 4.4}, we can restrict our consideration to Grassmannian, hyperquadric and spinor Grassmannian. If $\mathcal{S}_0$ is also a spinor Grassmannian, then the pair is of sub-diagram type, so admissible. 

Let's consider the case where $\mathcal{S}_0$ is a Grassmannian, by virtue of Prop \ref{Prop 4.3}, we should find out all the embedding $\psi:\mathbb{P}^{r-1}\times\mathbb{P}^{s-1}\hookrightarrow G(2,m-2)$\ \ such that restriction of Plucker embedding of $G(2,m-2)$\ \ into $\mathbb{P}^{m(m-1)/2-1}$\ \ to the image of $\psi$\ \ yields the Segre embedding of $\mathbb{P}^{r-1}\times\mathbb{P}^{s-1}$\ \ into $\mathbb{P}^{rs-1}$.\ \ If we express $\psi$\ \ in terms of Harish-Chandra coordinates, then 

\begin{center}
$\psi:(w_1,...,w_{r-1})\bullet (z_1,...,z_{s-1})\hookrightarrow M_{2,m-2}$
\end{center}

where $M_{2,m-2}$\ \ denotes a $2\times (m-2)$\ \ complex matrix, such that all the $w$-coordinate($z$-coordinate) are put in the first (second) row and no $w$-coordinate is put in the same column as $z$-coordinate. Thus we must have $(r-1)+(s-1)\leq m-2$\ \ i.e., $r+s\leq m$. 

On the other hand, if $\mathcal{S}_0$\ \ is some hyperquadric $Q^n$. we have the reduction process 

\begin{center}
$(Q^n,G^{\Rmnum{2}}(m,m))\rightarrow (Q^{n-2},G(2,m-2))\rightarrow (Q^{n-4},\mathbb{P}^1\times \mathbb{P}^{m-3})$
\end{center}

So we should determine all the embedding $\psi: Q^r\hookrightarrow \mathbb{P}^1\times\mathbb{P}^{s},\ \ s\geq 2$\ \ such that restriction of Segre embedding to the image of $\psi$\ \ yields the usual embedding of $Q^r$\ \ into $\mathbb{P}^{r+1}$.\ \ Actually, the only candidate is $Q^2\cong \mathbb{P}^1\times\mathbb{P}^1$\ \ regardless of what $s$\ \ is and $\psi$\ \ is induced by respective inclusion $\mathbb{P}^1\subset\mathbb{P}^1,\ \ \mathbb{P}^1\subset\mathbb{P}^s$.\ \ Again by virtue of Prop \ref{Prop 4.3}, we reverse the reduction process and get the admissible pair $(Q^6,G^{\Rmnum{2}}(m,m))$, which is of sub-diagram type. Note that $(Q^n,Q^6)\ \ n=2,3,4,5$\ \ are all admissible pairs, so are $(Q^n,G^{\Rmnum{2}}(m,m))\ \ n=2,3,4,5$ by transitivity of admissibility which is immediate from the definition of admissible pair.

\begin{flushright}
QED
\end{flushright}

The remaining task towards completing the classification is to determine all the admissible pairs $(\mathcal{S}_0,\Rmnum{5})$ and $(\mathcal{S}_0,\Rmnum{6})$. Again 3 of Lemma \ref{Lemma 4.2} and Lemma \ref{Lemma 4.4} allow us to only consider the candidates Grassmannian, hyperquadric and spinor Grassmannian. 

The reduction process for $\Rmnum{5}$: 
\begin{align*}
&\mathcal{S}_0=G(r,s),\ \ (G(r,s),\Rmnum{5})\rightarrow (\mathbb{P}^{r-1}\times\mathbb{P}^{s-1},G^{\Rmnum{2}}(5,5))\\
&\mathcal{S}_0=Q^n,\ \ (Q^n,\Rmnum{5})\rightarrow (Q^{n-2},G^{\Rmnum{2}}(5,5))\rightarrow (Q^{n-4},G(2,3))\rightarrow (Q^{n-6},\mathbb{P}^{1}\times\mathbb{P}^{2})\\
&\mathcal{S}_0=G^{\Rmnum{2}}(n,n),\ \ (G^{\Rmnum{2}}(n,n),\Rmnum{5})\rightarrow (G(2,n-2),G^{\Rmnum{2}}(5,5))\rightarrow (\mathbb{P}^{1}\times\mathbb{P}^{n-3},G(2,3))\ \ n\geq 5
\end{align*}

and the reduction process for $\Rmnum{6}$:
\begin{align*}
& \mathcal{S}_0=G(r,s),\ \ (G(r,s),\Rmnum{6})\rightarrow (\mathbb{P}^{r-1}\times\mathbb{P}^{s-1},\Rmnum{5})\\
& \mathcal{S}_0=Q^m,\ \ (Q^m,\Rmnum{6})\rightarrow (Q^{m-2},\Rmnum{5})\rightarrow\cdots\rightarrow (Q^{m-8},\mathbb{P}^{1}\times\mathbb{P}^{2})\\
& \mathcal{S}_0=G^{\Rmnum{2}}(m,m),\ \ (G^{\Rmnum{2}}(m,m),\Rmnum{6})\rightarrow (G(2,m-2),\Rmnum{5})\rightarrow (\mathbb{P}^{1}\times\mathbb{P}^{m-3},G^{\Rmnum{2}}(5,5))
\end{align*}

With the above reduction processes, we only need to establish the following three propositions to complete our classification:

\begin{proposition}
\label{Prop 4.6}
For $(Q^n,\mathbb{P}^1\times\mathbb{P}^2)$\ \ to be admissible pair of VMRT, $n=1,2$.
\end{proposition}

\begin{proposition}
\label{Prop 4.7} 
For $(\mathbb{P}^r\times\mathbb{P}^s,\Rmnum{5})$\ \ to be admissible pair of VMRT, (1) either $r=1,\ \ 1\leq s\leq 5$\ \ (2) or $r=2,\ \ s=2$.
\end{proposition}

\begin{proposition}
\label{Prop 4.8} 
For $(\mathbb{P}^r\times\mathbb{P}^s,G^{\Rmnum{2}}(5,5))$\ \ to be admissible pair of VMRT, $r=1,\ \ 1\leq s\leq 3$.
\end{proposition}

Assuming the above three propositions and note that we have already determined all the admissible pair $(\mathbb{P}^r\times\mathbb{P}^s,G(2,n))$ of VMRT in the argument of proof of Lemma \ref{Lemma 4.5}, we immediately have the following conclusions by means of reversing the reduction processes thanks to Prop \ref{Prop 4.3}:

\begin{lemma}
\label{Lemma 4.10} 

If $(\mathcal{S}_0,\Rmnum{5})$\ \ is an admissible pair, then:
\begin{enumerate}
\item $\mathcal{S}_0=G(2,s),\ \ 2\leq s\leq 4$ (the pair is of deletion type in category 2 when $s=2,3$ and of special type in category 3 when $n=4$);
\item $\mathcal{S}_0=Q^n,\ \ 2\leq n\leq 8$ (the pair is of sub-diagram type in category 1 when $n=8$ and is in category 4 when $n=2,3,...,7$);
\item $\mathcal{S}_0=G^{\Rmnum{2}}(5,5)$ (this pair is of deletion type in category 2).
\end{enumerate} 

\end{lemma}

\begin{lemma}
\label{Lemma 4.11} 
If $(\mathcal{S}_0,\Rmnum{6})$\ \ is an admissible pair, then:
\begin{enumerate}
\item $\mathcal{S}_0=G(2,s),\ \ 2\leq s\leq 6$ (the pair is of deletion type in category 2 when $n=2,3$, is of special type in category 3 when $n=5,6$, and is in category 4 when $n=4$);
\item $\mathcal{S}_0=G(3,3)$ (this pair is of special type in category 3);
\item $\mathcal{S}_0=Q^n,\ \ 2\leq n\leq 10$ (the pair is of sub-diagram type when $n=10$, and is in category 4 when $2\leq n\leq 9$);
\item $\mathcal{S}_0=G^{\Rmnum{2}}(n,n),\ \ n=5,6$ (the pair is of deletion type in category 2 when $n=5$, and is of special type in categoty 3 when $n=6$).
\end{enumerate}

\end{lemma}

Now we prove Prop \ref{Prop 4.6}-Prop \ref{Prop 4.8}

\textbf{Proof of Prop \ref{Prop 4.6}} Note that $Q^2\cong\mathbb{P}^1\times\mathbb{P}^1$,\ \ so $\psi:Q^2\hookrightarrow\mathbb{P}^1\times\mathbb{P}^2$\ \ can be induced by inclusion $\mathbb{P}^1\subset\mathbb{P}^1,\ \ \mathbb{P}^1\subset\mathbb{P}^2$.\ \ There exists no $\psi:Q^3\rightarrow\mathbb{P}^1\times\mathbb{P}^2$, simply because $Q^3$ is irreducible while $\mathbb{P}^1\times\mathbb{P}^2$ is reducible.

\begin{flushright}
QED
\end{flushright}

\textbf{Proof of Prop \ref{Prop 4.7}} The marked Dynkin diagram for $\Rmnum{5},\ \ \Rmnum{6}$\ \ are shown below respectively:

\vspace{10mm}

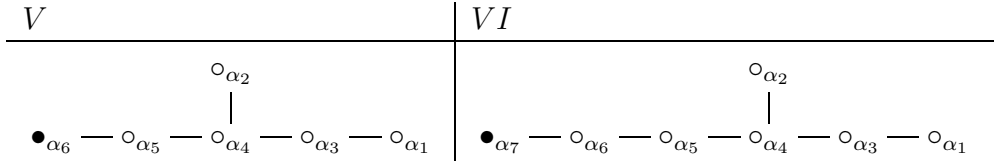
\begin{tabular}{l|l}
$\Rmnum{5}$ & $\Rmnum{6}$\\

\hline

\xymatrix@C=1em@R=1em{
&&\circ_{\alpha_2}\ar@{-}[d] \\
\bullet_{\alpha_6} \ar@{-}[r]&\circ_{\alpha_5}\ar@{-}[r]&\circ_{\alpha_4}\ar@{-}[r]&\circ_{\alpha_3}\ar@{-}[r]&\circ_{\alpha_1}
}

&

\xymatrix@C=1em@R=1em{
&&&\circ_{\alpha_2}\ar@{-}[d] \\
\bullet_{\alpha_7} \ar@{-}[r]&\circ_{\alpha_6}\ar@{-}[r]&\circ_{\alpha_5}\ar@{-}[r]&\circ_{\alpha_4}\ar@{-}[r]&\circ_{\alpha_3}\ar@{-}[r]&\circ_{\alpha_1}
}

\end{tabular}

\vspace{10mm}

It is well known from the geometry of cHSS that there are $16(27)$\ \ positive roots for $\Rmnum{5}(\Rmnum{6})$.\ \ Denote by $\beta_j,\ \ j=1,2,...,16$\ \ the $16$\ \ positive roots for $\Rmnum{5}$,\ \ then $16$\ \ of the $27$\ \ positive roots of $\Rmnum{6}$\ \ are obtained by attaching $\alpha_7$\ \ to $\beta_j$,\ \ i.e., $\gamma_j=\alpha_7+\beta_j,\ \ j=1,2,...,16$.

We start with $r=1$, we aim at showing $s\leq 5$, i.e., for any chosen positive roots $\gamma_{j_0}$ there exist at most 5 positive root $\gamma_{j_1},...,\gamma_{j_5}$ such that $\gamma_{j_k}+\gamma_{j_0}-\alpha_7,\ \ k=1,2,...,5$ remain to be positive roots, or equivalently, $\gamma_{j_k}-\gamma_{j_0}=\beta_{j_k}-\beta_{j_0},\ \ k=1,2,...,5$ are not roots. Actually this has been established by an early paper of Zhong Jia-Qing, cf \cite{[Zhong 84]}. Furthermore, the equality holds, i.e., there exists exactly five positive roots satisfying the aforementioned property when any $\gamma_{j_0}$ is chosen. For example, when $\beta_1=\alpha_6$, following the notation in [Zhong 84], we have $\beta_1^{\bot}=\{\alpha_6+2\alpha_5+2\alpha_4+\alpha_3+\alpha_2,\ \ \alpha_6+2\alpha_5+2\alpha_4+\alpha_3+\alpha_2+\alpha_1,\ \ \alpha_6+2\alpha_5+2\alpha_4+2\alpha_3+\alpha_2+\alpha_1,\ \ \alpha_6+2\alpha_5+3\alpha_4+2\alpha_3+\alpha_2+\alpha_1,\ \ \alpha_6+2\alpha_5+3\alpha_4+2\alpha_3+2\alpha_2+\alpha_1\}$. By transitivity, this is true for any other single root $\beta_j$.

Next we consider $r=2$ and $s\geq 2$ (note that $r=2,s=1$ is essentially the same as $r=1,s=2$ which has been established in the above argument). We aim at showing that $s=2$ is the only possibility, i.e., given any two roots $\beta_{j_1},\beta_{j_2}$ such that $\beta_{j_1}-\beta_{j_2}$ is a root (actually this is a compact root) or $\beta_{j_2}\notin \beta_{j_1}^{\bot}$ following the notation of \cite{[Zhong 84]}, then we have $\#(\beta_{j_1}^{\bot}\cap\beta_{j_2}^{\bot})=2$. This is essentially a checking game, we choose $\beta_{15}=\alpha_6+2\alpha_5+3\alpha_4+2\alpha_3+\alpha_2+\alpha_1,\ \ \beta_{16}=\alpha_6+2\alpha_5+3\alpha_4+2\alpha_3+2\alpha_2+\alpha_1$,\ \ one-by-one checking shows that $\beta_{15}^{\bot}\cap\beta_{16}^{\bot}=\{\alpha_6,\alpha_6+\alpha_5\}$.\ \ It remains to justify that this is true for all such pairs $(\beta_{j_1},\beta_{j_2})$\ \ if it holds for the above fixed pair. In fact, any such pair $(\beta_{j_1},\beta_{j_2})$\ \ determines a 2-dimensional subspace of tangent space $W\subset T_o(\Rmnum{5})$\  \ at reference point $o\in\Rmnum{5}$\ \ such that for any $v\in W$,\ \ the line $L$\ \ with $T_oL=v$\ \ shooting out from $o$\ \ is a minimal rational curve of $\Rmnum{5}$\ \ in some Harish-Chandra coordinate. So $W$\ \ determines a minimal rational curve on $\mathscr{C}_o(\Rmnum{5})\cong G^{\Rmnum{2}}(5,5)$. On the other hand, Aut$(G^{\Rmnum{2}}(5,5))$ acts transitively on the space of minimal rational curves of $G^{\Rmnum{2}}(5,5)$.

Then all the cases where $r,s\geq 3$ are impossible which follows immediately from the fact that $r=2, s\geq 3$ is impossible.

\begin{flushright}
QED
\end{flushright}

\textbf{Proof of Prop \ref{Prop 4.8}} Below is the marked Dynkin diagram of $G^{\Rmnum{2}}(5,5)$:
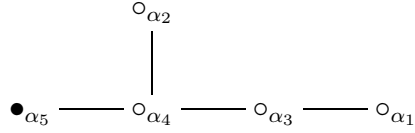
\begin{align*}
&
\xymatrix{
&\circ_{\alpha_2}\ar@{-}[d] \\
\bullet_{\alpha_5} \ar@{-}[r]&\circ_{\alpha_4}\ar@{-}[r]&\circ_{\alpha_3}\ar@{-}[r]&\circ_{\alpha_1}
}
&
\end{align*}

Denote by $\delta_k,\ \ k=1,2,...,10$ the 10 positive roots of $G^{\Rmnum{2}}(5,5)$, we need to show: 1. for any positive root $\delta_{k_0}$,\ \ Card$(\delta_{k_0}^{\bot})=3$; 2. for any two positive roots $\delta_{k_1},\delta_{k_2}$\ \ such that $\delta_{k_1}\notin\delta_{k_2}^{\bot},\ \ \delta_{k_1}^{\bot}\cap\delta_{k_2}^{\bot}=\emptyset$. The verification of the above two facts just borrows the argument for Prop \ref{Prop 4.7}.

\begin{flushright}
QED
\end{flushright}


\subsection{\textbf{Construction of standard embedding $i:\mathcal{S}_0\hookrightarrow\mathcal{S}$ for candidates of admissible pair $(\mathcal{S}_0,\mathcal{S})$}}
\label{section 4.2}

In the previous section, we theoretically list all the possible candidates of admissible pairs $(\mathcal{S}_0,\mathcal{S})$. We should explicitly  construct standard embedding $i:\mathcal{S}_0\hookrightarrow\mathcal{S}$ to show that the candidates are actually eligible to be admissible pairs. 

The strategy of our construction is a dichomy; for most cases, by virtue of the backgrounds provided in section \ref{intro}, in particular, Lemma \ref{Lemma 3 for paper 2}, Prop \ref{Prop 2 for paper 2}, we only need to construct a root correspondence $\Phi:\Delta_0\rightarrow\Delta$ such that; for the only three cases which are the item (a) (b) in category 3 in the Main Theorem 1, noting that all of them involve no exceptional space, we use coordinate to explicitly construct standard embedding. 

\medskip
	
	\textbf{Construction for category 1}
	
\medskip
	Pairs of sub-diagram type are naturally admissible pairs as was mentioned in Theorem \ref{Thm 1.1}.
	
\medskip

	\textbf{Construction for category 2}


	\begin{proposition}
	\label{Prop 4.12}
	The pair $(\mathcal{S}_0,\mathcal{S})$ of deletion type is an admissible pair.
	\end{proposition}
	
	\textbf{Proof} Suppose the marked Dynkin diagram $(\mathcal{D}(\mathcal{S}_0),\gamma_0)$\ \ is as follows: 


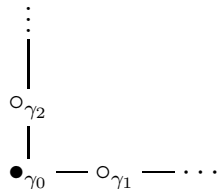
\centerline{\xymatrix@C=1em@R=1em{
\vdots\ar@{-}[d]\\
\circ_{\gamma_2}\ar@{-}[d] \\
\bullet_{\gamma_0} \ar@{-}[r]&\circ_{\gamma_1}\ar@{-}[r]&\cdots
}}

Note that since $\mathcal{S}_0$\ \ is a Hermitian Symmetric Space, then there are at most two nodes adjacent to the marked node $\gamma_0$.\ \ Then according to the definition, the marked Dynkin diagram $(\mathcal{D}(\mathcal{S}),\gamma)$\ \ should be like:

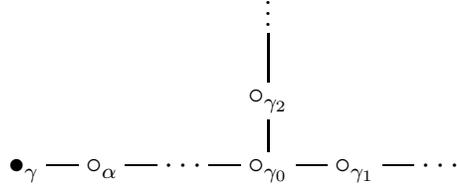
\centerline{\xymatrix@C=1em@R=1em{
&&&\vdots\ar@{-}[d]\\
&&&\circ_{\gamma_2}\ar@{-}[d] \\
\bullet_{\gamma}\ar@{-}[r]&\circ_{\alpha}\ar@{-}[r]&\cdots\ar@{-}[r]&\circ_{\gamma_0} \ar@{-}[r]&\circ_{\gamma_1}\ar@{-}[r]&\cdots
}}and we have an identification of $\mathcal{D}(\mathcal{S}_0)$\ \ with a sub-diagram of $\mathcal{D}(\mathcal{S})$\ \ not containing the node $\alpha$.\ \ We construct a mapping $\Phi:\mathcal{D}(\mathcal{S}_0)\rightarrow\Delta$ defined on simple roots as follows: $\Phi(\gamma_0)=\gamma$,  $\Phi(\gamma_1)=\alpha+\cdots+\gamma_0+\gamma_1$, $\Phi(\gamma_2)=\alpha+\cdots+\gamma_0+\gamma_2$, $\Phi(\beta)=\beta$, 
where $\beta$s are the rest nodes of $\mathcal{D}(\mathcal{S}_0)$\ \ which is identified with nodes of $\mathcal{D}(\mathcal{S})$\ \ in the sense of definition of deletion type pair. 

We can show that $\Phi$ preserves Cartan integers among simple roots, so it extends uniquely to isomorphism between root systems. Thanks to Proposition \ref{Prop 2 for paper 2}, we acquire the standard embedding $i_{\Phi}:\mathcal{S}_0\hookrightarrow\mathcal{S}$ induced by $\Phi$. 

\begin{flushright}
QED
\end{flushright}

\medskip

\textbf{Construction for category 3}

\medskip

\begin{enumerate}
\item $(G(4,2),\Rmnum{5})$

The marked Dynkin diagram of $G(4,2)$\ \ and $\Rmnum{5}$:
	
	\vspace{3mm}
	
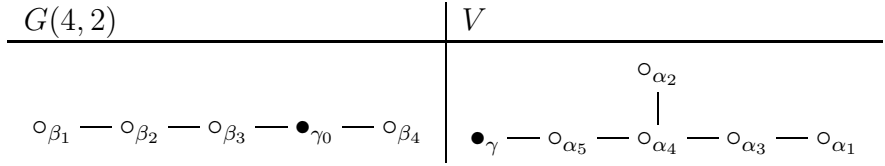
	\begin{tabular}{l|l}
	$G(4,2)$ & $\Rmnum{5}$ \\
	\hline
			\xymatrix@C=1em@R=1em{&&&& \\
	\circ_{\beta_1}\ar@{-}[r]&\circ_{\beta_2}\ar@{-}[r]&\circ_{\beta_3}\ar@{-}[r]&\bullet_{\gamma_0}\ar@{-}[r]&\circ_{\beta_4}
	}
	&
	\xymatrix@C=1em@R=1em{
&&\circ_{\alpha_2}\ar@{-}[d] \\
\bullet_{\gamma} \ar@{-}[r]&\circ_{\alpha_5}\ar@{-}[r]&\circ_{\alpha_4}\ar@{-}[r]&\circ_{\alpha_3}\ar@{-}[r]&\circ_{\alpha_1}
}
	
	\end{tabular}
	
	\vspace{5mm}
	
	Define $\Phi$ as:
	
	\vspace{3mm}
	
	 $\Phi(\gamma_0) =\gamma$, $\Phi(\beta_4)  = \alpha_5$, $\Phi(\beta_3)  = \alpha_5+2\alpha_4+\alpha_3+\alpha_2$, $\Phi(\beta_2)  = \alpha_1$, $\Phi(\beta_1)  = \alpha_3$.
	\vspace{8mm}

\item $(G^{\Rmnum{2}}(6,6),\Rmnum{6})$

The marked Dynkin diagram of $G^{\Rmnum{2}}(6,6)$ and $\Rmnum{6}$:
	
\vspace{3mm}

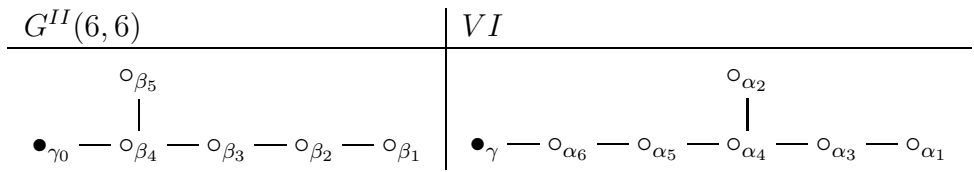
\begin{tabular}{l|l}
	$G^{\Rmnum{2}}(6,6)$ & $\Rmnum{6}$ \\
	\hline
	\xymatrix@C=1em@R=1em{
&\circ_{\beta_5}\ar@{-}[d] \\
\bullet_{\gamma_0} \ar@{-}[r]&\circ_{\beta_4}\ar@{-}[r]&\circ_{\beta_3}\ar@{-}[r]&\circ_{\beta_2}\ar@{-}[r]&\circ_{\beta_1}
}
	&
	
	\xymatrix@C=1em@R=1em{
&&&\circ_{\alpha_2}\ar@{-}[d] \\
\bullet_{\gamma} \ar@{-}[r]&\circ_{\alpha_6}\ar@{-}[r]&\circ_{\alpha_5}\ar@{-}[r]&\circ_{\alpha_4}\ar@{-}[r]&\circ_{\alpha_3}\ar@{-}[r]&\circ_{\alpha_1}
}
	
	\end{tabular}
	
	\vspace{5mm}
	
	Define $\Phi$ as: 
	
	\vspace{3mm}
	
	$\Phi(\gamma_0) = \gamma$, $\Phi(\beta_4)  = \alpha_6$, $\Phi(\beta_5)  = \alpha_5$, $\Phi(\beta_3)  = \alpha_5+2\alpha_4+\alpha_3+\alpha_2$, $\Phi(\beta_2)  = \alpha_1$, $\Phi(\beta_1)  = \alpha_3$.
	
		\vspace{8mm}

\item $(G(5,2),\Rmnum{6})$

The marked Dynkin diagram of $G(6,2)$ and $\Rmnum{6}$:
	
\vspace{3mm}

\begin{tabular}{l|l}
	$G(6,2)$ & $\Rmnum{6}$ \\
	\hline
			\xymatrix@C=1em@R=1em{&&&& \\
	\circ_{\beta_1}\ar@{-}[r]&\circ_{\beta_2}\ar@{-}[r]&\circ_{\beta_3}\ar@{-}[r]&\circ_{\beta_4}\ar@{-}[r]&\circ_{\beta_5}\ar@{-}[r]&\bullet_{\gamma_0}\ar@{-}[r]&\circ_{\beta_6}
	}
	&
	
	\xymatrix@C=1em@R=1em{
&&&\circ_{\alpha_2}\ar@{-}[d] \\
\bullet_{\gamma} \ar@{-}[r]&\circ_{\alpha_6}\ar@{-}[r]&\circ_{\alpha_5}\ar@{-}[r]&\circ_{\alpha_4}\ar@{-}[r]&\circ_{\alpha_3}\ar@{-}[r]&\circ_{\alpha_1}
}
	
	\end{tabular}
	
	\vspace{5mm}
	
	Define $\Phi$ as:
	
	\vspace{3mm}
	
	$\Phi(\gamma_0) = \gamma$, $\Phi(\beta_6)  = \alpha_6$, $\Phi(\beta_5)  = \alpha_6+2\alpha_5+2\alpha_4+\alpha_3+\alpha_2$, $\Phi(\beta_4)  = \alpha_1$, $\Phi(\beta_3)  = \alpha_3$, $\Phi(\beta_2)  = \alpha_4$, $\Phi(\beta_1)  = \alpha_2$.

	\vspace{8mm}

\item $(G(5,2),\Rmnum{6})$

$\Phi$ is obtained from the case of $(G(6,2),\Rmnum{6})$ by discarding $\Phi(\beta_1)=\alpha_2$.

	\vspace{8mm}
	
\item $(G(3,3),\Rmnum{6})$

The marked Dynkin diagram of $G(3,3)$ and $\Rmnum{6}$:
	
\vspace{3mm}

\begin{tabular}{l|l}
	$G(3,3)$ & $\Rmnum{6}$ \\
	\hline
			\xymatrix@C=1em@R=1em{&&&& \\
	\circ_{\beta_1}\ar@{-}[r]&\circ_{\beta_2}\ar@{-}[r]&\bullet_{\gamma_0}\ar@{-}[r]&\circ_{\beta_3}\ar@{-}[r]&\circ_{\beta_4}
	}
	&
	
	\xymatrix@C=1em@R=1em{
&&&\circ_{\alpha_2}\ar@{-}[d] \\
\bullet_{\gamma} \ar@{-}[r]&\circ_{\alpha_6}\ar@{-}[r]&\circ_{\alpha_5}\ar@{-}[r]&\circ_{\alpha_4}\ar@{-}[r]&\circ_{\alpha_3}\ar@{-}[r]&\circ_{\alpha_1}
}
	
	\end{tabular}
	
		\vspace{5mm}

	Define $\Phi$ as:

	\vspace{3mm}

	$\Phi(\gamma_0) = \gamma$, $\Phi(\beta_4)  = \alpha_4+\alpha_3$, $\Phi(\beta_3)  = \alpha_6+\alpha_5+\alpha_4+\alpha_3+\alpha_2+\alpha_1$, $\Phi(\beta_2)  = \alpha_6+\alpha_5+\alpha_4$, $\Phi(\beta_1)  = \alpha_5+\alpha_4+\alpha_3+\alpha_2$.
	
	\vspace{8mm}

It is not hard to check that all the above $\Phi$s are isomorphism between root systems such that $\Phi(\gamma_0)=\gamma$, thanks to Proposition \ref{Prop 2 for paper 2}, we know that they are all admissible pairs. We go on to deal with the rest two pairs (a) and (b) in category 3 to complete the construction:

	\vspace{8mm}
	
\item $(G^{\Rmnum{3}}(n,n),G(r,s)),\ \ 3\leq n\leq$min$\{r,s\}$

A Harish-Chandra coordinate chart on $G^{\Rmnum{3}}(n,n),G(r,s)$ can be represented by $n\times n$ symmetric matrix $\textbf{M}^{\Rmnum{3}}$ and $r\times s$ matrix $\textbf{M}$ respectively . Then the standard embedding $i:G^{\Rmnum{3}}(n,n)\hookrightarrow G(r,s)$ restricted to this Harish-Chandra coordinate is just canonically identifying $n\times n$ symmetric matrix with an upper-left $n\times n$-block in $r\times s$ matrix. The tangent mapping is the same in terms of matrix representation, i.e., 

\begin{center}
$M_{i,j}=\left\{\begin{array}{l}  \textbf{M}^{\Rmnum{3}}_{i,j},\ \ \textnormal{for}\ \ 1\leq i,j\leq n, \\ 0,\ \ \textnormal{for the rest entries}.  \end{array}\right.$
\end{center} 

Since $[\alpha]\in\mathscr{C}_o(G^{\Rmnum{3}}(n,n))$ if and only if $\alpha$ in matrix form is a rank 1 $n\times n$ symmetric matrix; $[\alpha]\in\mathscr{C}_o(G(r,s))$ if and only if $[\alpha]$ in matrix form is a rank 1 $r\times s$ matrix. Of course, such embedding is VMRT-respecting, for rank 1 symmetric matrix is mapped to rank 1  matrix and rank $> 1$ symmetric matrix is mapped to rank $>1$ matrix. 




\vspace{8mm}

\item $(G(r,s),G^{\Rmnum{2}}(n,n)),\ \ 3\leq r,s,\ \ r+s\leq n$

Harish-Chandra coordinate of $G^{\Rmnum{2}}(n,n)$ is represented by $n\times n$ complex skew-symmetric matrix $\textbf{M}^{\Rmnum{2}}$, the standard embedding $i:G(r,s)\hookrightarrow G^{\Rmnum{2}}(n,n)$ is essentially identifying $\{i,j\}$-th entry of $\textbf{M}$ with $\{i,j+r\}$-th entry of $\textbf{M}^{\Rmnum{2}}$  and skew-symmetrizing it. The tangent mapping is the same in terms of matrix representation, i.e., 

\begin{center}
$M^{\Rmnum{2}}_{i,j+r}=-M^{\Rmnum{2}}_{j+r,i}=\left\{\begin{array}{l}  \textbf{M}_{i,j},\ \ \textnormal{for}\ \ 1\leq i\leq r,1\leq j\leq s, \\ 0,\ \ \textnormal{for the rest entries}.  \end{array}\right.$
\end{center} 

Since $[\alpha]\in\mathscr{C}_o(G^{\Rmnum{2}}(n,n))$ if and only if $\alpha$ in matrix form is a rank 2 $n\times n$ skew-symmetric matrix; $[\alpha]\in\mathscr{C}_o(G(r,s))$ if and only if $[\alpha]$ in matrix form is a rank 1 $r\times s$ matrix. Of course, such embedding is VMRT-respecting, for rank 1 matrix is mapped to rank 2 skew-symmetric matrix and rank $> 1$ matrix is mapped to rank $>2$ skew-symmetric matrix.





\end{enumerate}

\subsection{\textbf{On the special pairs in category 3}}
\label{section 4.3}

In this section, we are interested in investigating the special pairs in category 3 in Main Theorem 1. Our result is:

\begin{theorem}
\label{Thm 4.13}
The special pairs $(\mathcal{S}_0,\mathcal{S})$ are all non-degenerate in the sense of Def \ref{Def 1.5}. As a result, if a locally closed submanifold $S\subset\mathcal{S}$ inherits a sub-VMRT structure modeled on $(\mathcal{S}_0,\mathcal{S})$, then $S\subset\mathcal{S}$ is line-preserving, i.e., for $\forall p\in S$, the minimal rational curve $L$ issuing from $p$ with $T_pL\in\widetilde{\mathscr{C}}_p(\mathcal{S}_0)=\widetilde{\mathscr{C}}_p(\mathcal{S})\cap T_p\mathcal{S}$ must lie on $S$. So $S\subset\mathcal{S}$ is an open subset of some algebraic submanifold when $\mathcal{S}\hookrightarrow\mathbb{P}^N$ by first canonical embedding.  
\end{theorem}

Note that, for these special pairs, we are able to at most answer the problem of algebraicity of the submanifolds $S\subset\mathcal{S}$ at present stage, instead of rigidity or non-rigidity. The major difficulty involved here is the failure of parallel transport of VMRT (cf. [HoM 10]) which is essential in establishing the rigidity. On the other hand, there is no clear clue to heuristically imply that they are non-rigid. So confirmation of rigidity/non-rigidity remains to be investigated in the future research. 

\medskip

We mainly focus on proving Thm \ref{Thm 4.13} in the rest of this section, the basis of our proof consist in the following well-known result in VMRT theory, the readers may refer to \cite{[HwM 01]}.

\begin{proposition}
\label{Prop 2.14 of paper 2}
Let $S\subset X$ be a submanifold of an RHS $X$. If any minimal rational curve of $X$ is also contained in $S$ whenever it is tangent to $S$ at some point, i.e., $S$ has the property of line-preserving, then $S$ is an algebraic.
\end{proposition}The property of line-preserving can deduced from the following result which is established in chapter 3 of \cite{[Zhang 14]},

\begin{proposition}
\label{Prop 2.15 of paper 2}
Let $(X_0,X)$ be an admissible pair of RHS, and a submanifold $S\subset X$ inherits a sub-VMRT structure modelled on $(X_0,X)$ (cf. Def \ref{Def 1.4}), if the pair $(X_0,X)$ is non-degenerate in the sense of Def \ref{Def 1.5}, then $S$ is line-preserving. 
\end{proposition}By virtue of the above two propositions, we only need to verify the non-degeneracy of the special type pairs $(\mathcal{S}_0,\mathcal{S})$. We first of all deal with $(G^{\Rmnum{3}}(n,n),G(r,s))$, $3\leq n\leq $min$\{r,s\}$. The embedding of respective abstract VMRT is 

\begin{align*}
\psi: & \mathbb{P}^{n-1}\hookrightarrow\mathbb{P}^{r-1}\times\mathbb{P}^{s-1},\\
      & (\zeta_1,...,\zeta_{n-1})\rightarrow (\zeta_1,...,\zeta_{n-1},0,...,0)\cdot (\zeta_1,...,\zeta_{n-1},0,...,0)
\end{align*} $\phi$ is Segre embedding

\begin{align*}
\phi: \mathbb{P}^{r-1}\times\mathbb{P}^{s-1} & \cong\mathscr{C}_o(G(r,s))\hookrightarrow\mathbb{P}^{rs-1}\\
(\xi_1,...,\xi_{r-1})\cdot (\eta_1,...,\eta_{s-1}) & \rightarrow (z_{i0}=\xi_i,z_{0j}=\eta_{j},z_{ij}=\xi_i\eta_j),\ \ 1\leq i\leq r-1,\ \ 1\leq j\leq s-1
\end{align*}It's not hard to compute the tangent space at the origin $\{0\}\in\mathbb{P}^{rs-1}$

\begin{center}
$T_{\{0\}}(\mathscr{C}_o(G(r,s)))=$span$\{\frac{\partial}{\partial z_{i0}},\frac{\partial}{\partial z_{0j}},\quad 1\leq i\leq r-1,\ \ 1\leq j\leq s-1\}$
\end{center} and $\phi_0=\phi\circ\psi$ is the embedding realized by $\Gamma(\mathbb{P}^{n-1},\mathcal{O}(2))$:

\begin{align*}
& \phi_0=\phi\circ\psi: \mathbb{P}^{n-1}\cong\mathscr{C}_o(G^{\Rmnum{3}}(n,n))\hookrightarrow\mathbb{P}^{(n+1)n/2-1}\subset\mathbb{P}^{rs-1}\\
& (\zeta_1,...,\zeta_{n-1})\rightarrow (z_{k0}=\zeta_k,z_{0l}=\zeta_{l},z_{kl}=\zeta_k\zeta_l),\ \ 1\leq k,l\leq n-1
\end{align*} Likewise, the tangent space at the origin

\begin{center}
$T_{\{0\}}(\mathscr{C}_o(G^{\Rmnum{3}}(n,n)))=$span$\{\frac{\partial}{\partial z_{0k}}+\frac{\partial}{\partial z_{k0}},\ \ 1\leq k\leq n-1\}$
\end{center}Denote by $\sigma_{\{0\}}$ the second fundamental form of $\mathscr{C}_o(G(r,s))$ at $\{0\}$, then it is not hard to show

\begin{align*}
\sigma_{\{0\}}(\frac{\partial}{\partial z_{i0}},\frac{\partial}{\partial z_{0j}}) & =\frac{\partial}{\partial z_{ij}},\ \ 1\leq i\leq r-1,\ \ 1\leq j\leq s-1;\\
\sigma_{\{0\}}(\frac{\partial}{\partial z_{i0}},\frac{\partial}{\partial z_{k0}}) & =\sigma_{\{0\}}(\frac{\partial}{\partial z_{0j}},\frac{\partial}{\partial z_{0l}})=0,\ \ 1\leq i,k\leq r-1,\ \ 1\leq j,l\leq s-1
\end{align*}Now let $v=\sum\limits_ia^i\frac{\partial}{\partial z_{i0}}+\sum\limits_jb^j\frac{\partial}{\partial z_{0j}}\in T_{\{0\}}(\mathscr{C}_o(G(r,s)))$ and the second fundamental form:

\begin{center}
$\sigma_{\{0\}}(v,\frac{\partial}{\partial z_{0k}}+\frac{\partial}{\partial z_{k0}})=(a^k+b^k)\frac{\partial}{\partial z_{kk}}+\sum\limits_{i\neq k}a^i\frac{\partial}{\partial z_{ik}}+\sum\limits_{j\neq k}b^j\frac{\partial}{\partial z_{kj}}$
\end{center} We just take $k=1,2$, and immediately get to the conclusion that for $\sigma_{\{0\}}$ to vanish identically for $k=1,2$ if and only if $v=0$. This shows nothing but the pair $(G^{\Rmnum{3}}(n,n),G(r,s))$ is non-degenerate. 

\medskip

Next we deal with $(G(r,s),G^{\Rmnum{2}}(n,n)),\ \ 3\leq r,s\ \ r+s\leq n$.

\bigskip The embedding of respective abstract VMRT $\psi:\mathbb{P}^{r-1}\times\mathbb{P}^{s-1}\hookrightarrow G(2,n-2)$ is

\begin{center}

$(\xi_1,...,\xi_{r-1})\cdot (\eta_1,...,\eta_{s-1})\rightarrow 
$$\left[\begin{array}{cccccccc}
\xi_1 & \cdots & \xi_{r-1} &  0     & \cdots &   0        &  0 & \cdots\\
   0  & \cdots &    0      & \eta_1 & \cdots & \eta_{s-1} &  0 & \cdots
\end{array}\right]$$ $
\end{center} $\phi$ is Pl$\ddot{\textnormal{u}}$cker embedding:

\begin{center}
$\phi:  G(2,n-2)\cong\mathscr{C}_o(G^{\Rmnum{2}}(n,n)) \hookrightarrow\mathbb{P}^{n(n-1)/2-1}$\\

\end{center}

\begin{center}
$ $$\left[\begin{array}{ccc}
z_{11} & \cdots & z_{1,n-2}\\
z_{21} & \cdots & z_{2,n-2}
\end{array}\right]$$ \rightarrow (z_{2i}\epsilon_1\wedge e_i, z_{1j}\epsilon_2\wedge e_j;\quad (z_{1i}z_{2j}-z_{1j}z_{2i})e_i\wedge e_j),\ \ 1\leq i,j\leq n-2$

\end{center} It's not hard to compute the tangent space at the origin $\{0\}\in\mathbb{P}^{n(n-1)/2-1}$

\begin{center}
$T_{\{0\}}(\mathscr{C}_o(G^{\Rmnum{2}}(n,n)))=$span$\{\epsilon_1\wedge e_i,\epsilon_2\wedge e_j,\ \ 1\leq i,j\leq n-2\}$
\end{center} and $\phi_0=\phi\circ\psi$ is the Segre embedding, likewise we can compute the tangent space of $\mathscr{C}_o(G(r,s))\subset\mathbb{P}^{rs-1}\subset\mathbb{P}^{n(n-1)/2-1}$ at origin:

\begin{center}
$T_{\{0\}}(\mathscr{C}_o(G(r,s)))=$span$\{\epsilon_1\wedge e_l,\epsilon_2\wedge e_k,\ \ 1\leq k\leq r-1, r\leq l\leq r+s-2\}$
\end{center} Concerning the second fundamental form of $\mathscr{C}_o(G^{\Rmnum{2}}(n,n))$ at origin, it is not hard to show

\begin{align*}
& \sigma_{\{0\}}(\epsilon_1\wedge e_i,\epsilon_2\wedge e_j)=-e_i\wedge e_j,\ \ 1\leq i<j\leq n-2\\
& \sigma_{\{0\}}(\epsilon_1\wedge e_i,\epsilon_1\wedge e_j)=\sigma_{\{0\}}(\epsilon_2\wedge e_i,\epsilon_2\wedge e_j)=0
\end{align*} Now choose some $v=\sum\limits_ia^i \epsilon_1\wedge e_i+\sum\limits_j b^j\epsilon_2\wedge e_j\in T_{\{0\}}(\mathscr{C}_o(G^{\Rmnum{2}}(n,n)))$ and 

\begin{align*}
& \sigma_{\{0\}}(v,\epsilon_1\wedge e_l)=\sum\limits_{j\neq l}b^j e_j\wedge e_l,\ \ r\leq l\leq r+s-2\\
& \sigma_{\{0\}}(v,\epsilon_2\wedge e_k)=\sum\limits_{j\neq l}a^i e_i\wedge e_k,\ \ 1\leq k\leq r-1
\end{align*} Since $r,s\geq 3$, there are more than one $l$s and $k$s as subindices, for the above second fundamental form to vanish for every $l,k$ if and only if $v=0$, this shows nothing but non-degeneracy of the pair $(G(r,s),G^{\Rmnum{2}}(n,n)),\ \ 3\leq r,s\ \ r+s\leq n$. 

\bigskip

We proceed to deal with rest 5 pairs in category 3. A common characteristic of them is that standard embedding of each of them $i_{\Phi}:\mathcal{S}_0\hookrightarrow\mathcal{S}$ is induced by root correspondence. Thus we can verify their non-degeneracy through the tedious one-by-one examination of the root vectors based on Lemma 3.5 in [HoM 10], we demonstrate the checking for the pair $(G(3,3),\Rmnum{6})$ and thus confirm its non-degeneracy as follows.

The standard embedding $i_{\Phi}:G(3,3)\hookrightarrow\Rmnum{6}$ induced by the $\Phi$ as defined in subsection \ref{section 4.2} gives rise to the embedding between VMRTs d$i_{\Phi}:\mathbb{P}^2\times\mathbb{P}^2\cong\mathscr{C}_o(G(3,3))\hookrightarrow\mathscr{C}_o(\Rmnum{6})\cong\Rmnum{5}$. Denote by $[\alpha]=[E_{\gamma}]\in\mathscr{C}_o(\Rmnum{6})$, then $\mathscr{H}_{\alpha}=T_{[\alpha]}(\mathscr{C}_o(\Rmnum{6}))=$span$\{u_1,...,u_{16}\}$, the normal complement $\mathscr{N}_{\alpha}=$span$\{v_1,...,v_{10}\}$ (cf. Proposition \ref{Prop 1.8 for paper 2}), the positive root vectors $u_i$s and $v_i$s are listed as below on the left and on the right respectively (cf. Lemma \ref{Lemma 2 for paper 2}):

\bigskip

	\begin{tabular}{l|l}
$u_1=E_{\gamma+\alpha_6}$,                                            & $v_1=E_{\gamma+2\alpha_6+2\alpha_5+2\alpha_4+\alpha_3+\alpha_2}$, \\
$u_2=E_{\gamma+\alpha_6+\alpha_5}$,                                   & $v_2=E_{\gamma+2\alpha_6+2\alpha_5+2\alpha_4+\alpha_3+\alpha_2+\alpha_1}$,\\
$u_3=E_{\gamma+\alpha_6+\alpha_5+\alpha_4}$,                        & $v_3=E_{\gamma+2\alpha_6+2\alpha_5+2\alpha_4+2\alpha_3+\alpha_2+\alpha_1}$,\\
$u_4=E_{\gamma+\alpha_6+\alpha_5+\alpha_4+\alpha_2}$,               & $v_4=E_{\gamma+2\alpha_6+2\alpha_5+3\alpha_4+2\alpha_3+\alpha_2+\alpha_1}$,\\
$u_5=E_{\gamma+\alpha_6+\alpha_5+\alpha_4+\alpha_3}$,             & $v_5=E_{\gamma+2\alpha_6+2\alpha_5+3\alpha_4+2\alpha_3+2\alpha_2+\alpha_1}$,\\
$u_6=E_{\gamma+\alpha_6+\alpha_5+\alpha_4+\alpha_3+\alpha_1}$,    & $v_6=E_{\gamma+2\alpha_6+3\alpha_5+3\alpha_4+2\alpha_3+\alpha_2+\alpha_1}$,\\
$u_7=E_{\gamma+\alpha_6+\alpha_5+\alpha_4+\alpha_3+\alpha_2}$,    & $v_7=E_{\gamma+2\alpha_6+3\alpha_5+3\alpha_4+2\alpha_3+2\alpha_2+\alpha_1}$,\\
$u_8=E_{\gamma+\alpha_6+\alpha_5+2\alpha_4+\alpha_3+\alpha_2}$,   & $v_8=E_{\gamma+2\alpha_6+3\alpha_5+4\alpha_4+2\alpha_3+2\alpha_2+\alpha_1}$,\\
$u_9=E_{\gamma+\alpha_6+2\alpha_5+2\alpha_4+\alpha_3+\alpha_2}$,  & $v_9=E_{\gamma+2\alpha_6+3\alpha_5+4\alpha_4+3\alpha_3+2\alpha_2+\alpha_1}$,\\
$u_{10}=E_{\gamma+\alpha_6+\alpha_5+\alpha_4+\alpha_3+\alpha_2+\alpha_1}$,  & $v_{10}=E_{\gamma+2\alpha_6+3\alpha_5+4\alpha_4+3\alpha_3+2\alpha_2+2\alpha_1}$,                 \\
$u_{11}=E_{\gamma+\alpha_6+\alpha_5+2\alpha_4+\alpha_3+\alpha_2+\alpha_1}$, &                 \\
$u_{12}=E_{\gamma+\alpha_6+\alpha_5+2\alpha_4+2\alpha_3+\alpha_2+\alpha_1}$,&                 \\
$u_{13}=E_{\gamma+\alpha_6+2\alpha_5+2\alpha_4+\alpha_3+\alpha_2+\alpha_1}$,  & \\
$u_{14}=E_{\gamma+\alpha_6+2\alpha_5+2\alpha_4+2\alpha_3+\alpha_2+\alpha_1}$, & \\
$u_{15}=E_{\gamma+\alpha_6+2\alpha_5+3\alpha_4+2\alpha_3+\alpha_2+\alpha_1}$, & \\
$u_{16}=E_{\gamma+\alpha_6+2\alpha_5+3\alpha_4+2\alpha_3+2\alpha_2+\alpha_1}$,& \\
	
	\end{tabular}
	
	\medskip

Via the embedding d$i_{\Phi}$, we identify $T_{[\alpha]}(G(3,3))$ with a subspace of  $ T_{[\alpha]}(\Rmnum{6})$, then $T_{[\alpha]}(G(3,3))$=span$\{u_3,u_9,u_{10},u_{12}\}$. Now we can check the non-degeneracy by virtue of Lemma 3.5 in [HoM 10]. (Note that we only need to check those root vectors in $T_{[\alpha]}(\Rmnum{6})\backslash T_{[\alpha]}(G(3,3))$) 

\medskip
\centerline{$\sigma(u_1,u_9)=v_1\neq 0$, $\sigma(u_2,u_{12})=v_3\neq 0$, $\sigma(u_4,u_{12})=v_5\neq 0$,}
\centerline{$\sigma(u_5,u_{10})=v_2\neq 0$, $\sigma(u_6,u_9)=v_6\neq 0$, $\sigma(u_7,u_3)=v_1\neq 0$,}
\centerline{$\sigma(u_8,u_{10})=v_5\neq 0$, $\sigma(u_{11},u_9)=v_8\neq 0$, $\sigma(u_{13},u_{12})=v_{10}\neq 0$,}
\centerline{$\sigma(u_{14},u_3)=v_6\neq 0$,$\sigma(u_{15},u_{10})=v_{10}\neq 0$, $\sigma(u_{16},u_3)=v_9\neq 0$}

\medskip
We can perform the above checking procedure for the rest pairs listed as (c) (d) (e) (g) in Main Theorem 1.1, but we would rather choose a more conceptual and geometric approach. For this we begin with the following lemma:


\begin{lemma}
\label{Lemma 4.14}
Let $(\mathcal{S}_0,\mathcal{S})$ be admissible pair of cHSS, $\widetilde{\mathscr{C}}_o(\mathcal{S}_0)=\widetilde{\mathscr{C}}_o(\mathcal{S})\cap T_o\mathcal{S}$. Let $\zeta(\xi)=\alpha+\xi+\sigma'(\xi,\xi)$ be vector valued quadratic polynomial defined on $\xi\in T_{\alpha}(\widetilde{\mathscr{C}}_o(\mathcal{S}))$ (cf. Lemma \ref{Lemma 1 of paper 2}). Suppose some $\eta\in T_{\alpha}(\widetilde{\mathscr{C}}_o(\mathcal{S}))$ such that $\sigma'(\eta,\eta)=0$ and $\sigma'(\eta,T_{\alpha}(\widetilde{\mathscr{C}}_o(\mathcal{S}_0)))=0$, then for any $\alpha'\in\widetilde{\mathscr{C}}_o(\mathcal{S}_0)$ and any $\lambda\in\mathbb{C},\ \ \alpha'+\lambda\eta\in\widetilde{\mathscr{C}}_o(\mathcal{S})$, i.e., the unique projective line $L\subset\mathbb{P}(T_o\mathcal{S})$ passing through $[\alpha']$ and $[\eta]$ also lies on $\mathscr{C}_o(\mathcal{S})$ for any $[\alpha']\in\mathscr{C}_o(\mathcal{S}_0)$ 
\end{lemma}

\textbf{Proof} Since sub-VMRT is obtained by linear intersection, then the restriction of $\zeta$ to $T_{\alpha}(\widetilde{\mathscr{C}}_o(\mathcal{S}_0))$ yields the quasi-projective submanifold whose Zariski closure in $\mathbb{P}(T_o\mathcal{S}_0)$ is $\mathscr{C}_o(\mathcal{S}_0)$. Then for $\alpha'\in\widetilde{\mathscr{C}}_o(\mathcal{S}_0),\ \ \exists\xi'\in T_{\alpha}(\widetilde{\mathscr{C}}_o(\mathcal{S}_0))$ s.t. $\zeta(\xi')=\alpha+\xi'+\sigma'(\xi',\xi')=\alpha'$. As a result, $\zeta(\xi'+\lambda\eta)=\alpha+(\xi'+\lambda\eta)+\sigma'(\xi'+\lambda\eta,\xi'+\lambda\eta)$. Using $\sigma'(\xi',\eta)=\sigma'(\eta,\eta)=0$, we get $\zeta(\xi'+\lambda\eta)=\alpha+\xi'+\sigma'(\xi',\xi')+\lambda\eta=\alpha'+\lambda\eta\in\widetilde{\mathscr{C}}_o(\mathcal{S})$.

\begin{flushright}
QED
\end{flushright}

If $\phi:\mathcal{A}\stackrel{\cong}{\rightarrow}\mathscr{C}_o(\mathcal{S})\subset\mathbb{P}(T_o\mathcal{S})$ is realized through first canonical embedding, then $\sigma'(\eta,\eta)=0\Rightarrow[\eta]\in$ ``the VMRT of $\mathscr{C}_o(\mathcal{S})$'', denoted by $\mathscr{C}^2_o(\mathcal{S})$. We have:

\begin{proposition}
\label{Prop 4.15}
Let $(\mathcal{S}_0,\mathcal{S})$ be admissible pair of cHSS such that both abstract VMRTs $\mathcal{A}_0,\mathcal{A}$ are embedded by first canonical embedding, if it is degenerate and the degeneracy is caused by $\eta$ which satisfies $\sigma'(\eta,\eta)=0$ as in Lemma \ref{Lemma 4.14}, then there exists a (global) holomorphic embedding $i_{*}:\mathscr{C}_o(\mathcal{S}_0)\hookrightarrow\mathscr{C}^2_o(\mathcal{S})$ such that minimal rational curves of $\mathscr{C}_o(\mathcal{S}_0)$ are are mapped to those of $\mathscr{C}^2_o(\mathcal{S})$.
\end{proposition}

\textbf{Proof} As implied by Lemma \ref{Lemma 4.14}, the projective line $L$ linking $[\alpha']$ and $[\eta]$ lies on $\mathscr{C}_o(\mathcal{S})$ as well for $\forall [\alpha']\in\mathscr{C}_o(\mathcal{S}_0)$, then its tangent vector $T_{[\eta]}L\in\mathscr{C}^2_o(\mathcal{S})$ is uniquely determined by $L$, thus a holomorphic embedding $i_*:\mathscr{C}_o(\mathcal{S}_0)\hookrightarrow\mathscr{C}^2_o(\mathcal{S})$ is established. Choose an arbitrary projective line $L'\subset\mathscr{C}_o(\mathcal{S}_0)$, and connects $[\alpha']$ with $[\eta]$ by the unique projective line $L_{[\alpha']}\subset\mathscr{C}_o(\mathcal{S})$ for $\forall[\alpha']\in L'$, thus some $\mathbb{P}^2$ is determined in $\mathbb{P}(T_o\mathcal{S})$, this $\mathbb{P}^2$ corresponds to a projective line in $\mathscr{C}^2_o(\mathcal{S})$. 

\begin{flushright}
QED
\end{flushright}

Note that for the five pairs labeled (c)-(g) in category 3 of the Main Theorem 1, standard embedding of each of them comes from a root correspondence $\Phi$ (Definition \ref{def 6 for paper 2}), if degeneracy occurs among them, then it must be caused by some root vector, i.e.,

\begin{lemma}
\label{Lemma 4.15'}
Let $(\mathcal{S}_0,\mathcal{S})$ be degenerate admissible pair associated to marked Dynkin diagrams $(\mathcal{D}(G_0),\gamma_0),(\mathcal{D}(G),\gamma)$, with the sets of positive roots $E_0^+=\{\alpha_0\}\sqcup H_{\alpha_0}\sqcup N_{\alpha_0},\ \ E^+=\{\alpha\}\sqcup H_{\alpha}\sqcup N_{\alpha}$, respectively. Suppose standard embedding $i_{\Phi}:\mathcal{S}_0\hookrightarrow\mathcal{S}$ is induced by root correspondence $\Phi:\Delta_0\rightarrow\Delta$, regard $E^+_0$ as a subset of $E^+$ via the homomorphism $\widetilde{\Phi}$ (cf. Lemma \ref{Lemma 3 for paper 2}), then there must exists some root vectors $E_{\gamma+\theta_{n+1}}\in H_{\alpha}\backslash H_{\alpha_0}$, such that $E_{\gamma+\theta_{n+1}}\in$Ker$'(\sigma')$ as in Def \ref{Def 1.5}.
\end{lemma}

\textbf{Proof} We have $\mathscr{H}_{\alpha_0}=$span$\{H_{\alpha_0}\},\ \ \mathscr{H}_{\alpha}=$span$\{H_{\alpha}\}$. We make use of Lemma 3.5 in [HoM 10], which implies that taking second fundamental form results in a shift among root vectors and the root vectors are mutually linearly independent, then this lemma follows immediately. 

\begin{flushright}
QED
\end{flushright}

In particular, let $(\mathcal{S}_0,\mathcal{S})$ be one of the speical type pairs (c)-(g) in category 3 of Main Theorem 1, as we have shown in subsection \ref{section 4.2}, the standard embedding $i_{\Phi}:\mathcal{S}_0\hookrightarrow\mathcal{S}$ is induced by root correspondence $\Phi$. So if $(\mathcal{S}_0,\mathcal{S})$ is degenerate, the degeneracy must be caused by some root vector $E_{\gamma+\theta_{n+1}}$. On the other hand, again in view of Lemma 3.5 in [HoM 10], we have $\sigma'(E_{\gamma+\theta_{n+1}},E_{\gamma+\theta_{n+1}})=0$ for $\gamma+2\theta_{n+1}$ is never a root. As a result, Proposition \ref{Prop 4.15} implies the existence of $i_*:\mathscr{C}_o(\mathcal{S}_0)\hookrightarrow\mathscr{C}^2_o(\mathcal{S})$. We aim at arguing by contradiction to show there exists no such embedding $i_*$ for the pairs $(\mathcal{S}_0,\mathcal{S})$, then non-degenracy follows. Specifically, we need to show:


\begin{proposition}
\label{Prop 4.16}
There exists no global holomorphic embedding $i_*:\mathcal{A}_0\hookrightarrow\mathcal{A}$ preserving minimal rational curve if $(\mathcal{A}_0=\mathbb{P}^1\times\mathbb{P}^3,\mathcal{A}=G(2,3))$ or $(\mathcal{A}_0=\mathbb{P}^1\times\mathbb{P}^4,\mathcal{A}=G^{\Rmnum{2}}(5,5))$.
\end{proposition}

Assuming Prop \ref{Prop 4.16}, we proceed to confirm the non-degeneracy of the five pairs (c)-(g):

\begin{enumerate}

\item If $(G(2,4),\Rmnum{5})$ is degenerate, then there exists $i_*:\mathbb{P}^1\times\mathbb{P}^3\hookrightarrow G(2,3)$, contradiction arises.

\item If $(G(2,5),\Rmnum{6})$ is degenerate, then there exists $i_*:\mathbb{P}^1\times\mathbb{P}^4\hookrightarrow G^{\Rmnum{2}}(5,5)$, contradiction arises.  

\item If $(G^{\Rmnum{2}}(6,6),\Rmnum{6})$ is degenerate, then there exists $i_*:G(4,2)\rightarrow G^{\Rmnum{2}}(5,5)$, which can be reduced further to $i_{**}:\mathbb{P}^1\times\mathbb{P}^3\hookrightarrow G(2,3)$ because $i_*$ preserves M.R.C it is VMRT preserving. Contradiction arises. 

\item By (2), since $G(2,5)\hookrightarrow G(2,6)$, we have $(G(2,6),\Rmnum{6})$ is also non-degenerate.

\end{enumerate}

\textbf{Proof of Prop \ref{Prop 4.16}} 
\begin{itemize}
\item
For $\mathcal{A}_0=\mathbb{P}^1\times\mathbb{P}^3,\mathcal{A}=G(2,3)$, select any fixed $x\in\mathbb{P}^3,\ \ f(\mathbb{P}^1,x)$ is a minimal rational curve, then $(t,x)\in\mathbb{P}^1\times\mathbb{P}^3$, $\frac{\partial}{\partial t}f(t,x)|_{t=0}$ defines a nowhere vanishing section $\eta$ of the normal bundle $N_{S|G(2,3)}$ where $S=f(0,\mathbb{P}^3)\subset G(2,3)$. On the other hand, $G(2,3)$ admits Grassmannian structure, i.e., $T_{G(2,3)}\cong U\otimes V$ where $U,V$ are universal vector bundles of dimension $2,3$ respectively. Since any line $l\in\mathbb{P}^3$ is mapped to minimal rational curve under $f(0,l)$, without loss of generality, we may assume that 

\begin{center}
$f(0;\mathbb{P}^3):(z_1,z_2,z_3)\hookrightarrow 
$$\left[\begin{array}{ccc}
0 & 0 & 0\\
z_1 & z_2 & z_3
\end{array}\right]$$, $
\end{center} then we quickly have $U|_S=\mathcal{O}\oplus\mathcal{O}(1),\ \ V|_S=T_{\mathbb{P}^3}(-1)$, so $T_{G(2,3)}|_S=(\mathcal{O}\oplus\mathcal{O}(1))\otimes T_{\mathbb{P}^3}(-1)=T_{\mathbb{P}^3}\oplus T_{\mathbb{P}^3}(-1)$, i.e., the normal bundle $N_{S|G(2,3)}=T_{\mathbb{P}^3}(-1)$.

Tensoring $\eta$ with some holomorphic section of $\mathcal{O}(1)$, we get a holomorphic section of $T_{\mathbb{P}^3}$ which only vanishes along a hyperplane to order 1, contradicting the fact that every section of $T_{\mathbb{P}^3}$ is an Euler vector field whose vanishing loci strictly contains a hyperplane.

\vspace{5mm}

\item
For $\mathcal{A}_0=\mathbb{P}^1\times\mathbb{P}^4,\mathcal{A}=G^{\Rmnum{2}}(5,5)$, in the same sense as the above case, we produce a nowhere vanishing section $\eta$ of the normal bundle $N_{S|G^{\Rmnum{2}}(5,5)}$ where $S$ is the image $f(0;\mathbb{P}^4)\subset G^{\Rmnum{2}}(5,5)$, we seek for contradiction arising from this nowhere vanishing section $\eta$. Consider the following map: 

\begin{center}
$\mathbb{P}^4\stackrel{d}{\hookrightarrow}\mathbb{P}^4\times\mathbb{P}^4\stackrel{i}{\hookrightarrow} G(5,5)$,
\end{center} where $d$ is diagonal embedding and $i$ puts one of the factors of the product in row and the other in column of the $5\times 5$-matrix, we denote the image by $\mathbb{P}^4_r,\mathbb{P}_c^4$ respectively, i.e.,

\begin{center}
$(z_1,z_2,z_3,z_4)\cdot (w_1,w_2,w_3,w_4)\hookrightarrow
$$\left[\begin{array}{ccccc}
0 & z_1 & z_2 & z_3 & z_4\\
w_1 & 0 & 0 & 0 & 0\\
w_2 & 0 & 0 & 0 & 0\\
w_3 & 0 & 0 & 0 & 0\\
w_4 & 0 & 0 & 0 & 0\\
\end{array}\right]$$. $
\end{center} Let $U$ be the $5$-dimensional universal bundle on $G(5,5)$, we have $U_{\mathbb{P}_r^4}=\mathcal{O}^4\oplus\mathcal{O}(1),\ \ U_{\mathbb{P}^4_c}=T_{\mathbb{P}^4}(-1)\oplus\mathcal{O}$. Note that $S\subset G^{\Rmnum{2}}(5,5)\subset G(5,5)$ can be regarded as the image $i\circ d(\mathbb{P}^4)\subset G(5,5)$. Through the pull-back $d^*$, we get $U_{S}=T_{\mathbb{P}^4}(-1)\oplus\mathcal{O}(1)$. Taking the wedge product of $U$ and restricting it to $S$, we have $T_{G^{\Rmnum{2}}(5,5)}|_S=T_{\mathbb{P}^4}\oplus \Lambda^2(T_{\mathbb{P}^4}(-1))$. So we get the normal bundle $N_{S|G^{\Rmnum{2}}(5,5)}=\Lambda^2(T_{\mathbb{P}^4}(-1))$. 

Another fact about the section $\eta$ is that for $\forall x\in S\subset G^{\Rmnum{2}}(5,5),\ \ \eta(x)$ is a ''decomposable`` vector, i.e., $\eta(x)=e_1\wedge e_2$ for linearly independent $e_1,e_2\in T_{\mathbb{P}^4}(-1)_x$. So the section $\eta$ gives rise a rank $2$ subbundle $W\subset T_{\mathbb{P}^4}(-1)$. Thus $T_{\mathbb{P}^4}(-1)$ is topologically isomorphic to the sum of two rank $2$ subbundles $W\oplus V$, so that we have by Cartan formula the total Chern class relation $c(T_{\mathbb{P}^4}(-1))=c(W)c(V)$. The total Chern class of $T_{\mathbb{P}^4}(-1)$ may be computed through Euler sequence:

\begin{center}
$0\rightarrow \mathcal{O}(-1)\rightarrow \mathcal{O}^{n+1}\rightarrow T_{\mathbb{P}^n}(-1)\rightarrow 0$,
\end{center} we have $c(T_{\mathbb{P}^4}(-1))=1+\delta+\delta^2+\delta^3+\delta^4$, where the $\delta$s denote the generators of one-dimensional cohomology groups of $\mathbb{P}^4$. Clearly, $c(T_{\mathbb{P}^4}(-1))$ does not admit a decomposition into two polynomials about $\delta$ with integral coefficients, we get to contradiction. 

\end{itemize}

\begin{flushright}
QED
\end{flushright}

\textbf{Remark} In a similar fashion, we reasonably expect to establish the non-existence of holomorphic embedding $f:\mathbb{P}^2\times\mathbb{P}^2\hookrightarrow G^{\Rmnum{2}}(5,5)$ preserving minimal rational curves, so that we can also adopt a geometric and conceptual approach to confirming the non-degeneracy of the pair $(G(3,3),\Rmnum{6})$ without having to check root vectors one-by-one as we have done before (note that $\mathscr{C}^2_o(\Rmnum{6})\cong G^{\Rmnum{2}}(5,5)$). A source of relevant knowledge about vector bundles on projective spaces may be \cite{[OSS 80]}


\subsection{\textbf{Phenomenon of non-rigidity}}
\label{section 4.4}

Pairs $(Q^n,Q^m)$ of hyperquadrics provide typical examples of non-rigid pairs, since every submanifold $S\subset Q^m$ must inherit a sub-VMRT structure modeled on $(Q^n,Q^m)$, note that $\mathbb{P}(V)\cap\mathscr{C}_o(Q^m)\cong Q^{\textnormal{dim}(V)-2}$ for generic projective subspace $\mathbb{P}(V)\subset\mathbb{P}(T_o(Q^m))$. This section is designed for proving Main Theorem 2, as a sufficient condition for pairs to be non-rigid. The following preparatory lemma is implicitly contained in the procedure of classifying admissible pairs of cHSS: 

\begin{lemma}
\label{Prop 4.17}
Let $(\mathcal{S}_0,\mathcal{S})$ be a degenerate admissible pair of cHSS where $\mathcal{S}_0$ is allowed to be linear, then except the three pairs $(Z^{n-1}_{\textnormal{max}},Q^{2n-1}),\ \ n\geq 2$ and $(Z^1_{\textnormal{max}},G^{\Rmnum{3}}(n,n)),\ \ n\geq 2$ and $(Q^n,Q^m),\ \ n\equiv m$(mod 1), the standard embedding $i_{\Phi}:\mathcal{S}_0\hookrightarrow\mathcal{S}$ must be induced by root correspondence $\Phi:\Delta_0\rightarrow\Delta$. 
\end{lemma}

\textbf{Proof} If $\mathcal{S}_0$ is not linear, then from our classification theorem established as Main Theorem 1, the only admissible pairs which cannot be induced by root correspondence are $(Q^n,Q^m),\ \ n\equiv m$(mod 1) and $(G^{\Rmnum{3}}(n,n),G(r,s)),\ \ 3\leq n\leq$min$\{r,s\}$. Note that the latter pair is non-degenerate as we have verified in subsection \ref{section 4.3}. If $\mathcal{S}_0$ is linear, then from Choe and Hong 's result ([HoC 04]), $(Z^{n-1}_{\textnormal{max}},Q^{2n-1}),\ \ n\geq 2$ and $(Z^1_{\textnormal{max}},G^{\Rmnum{3}}(n,n)),\ \ n\geq 2$ are degenerate. 

\begin{flushright}
QED
\end{flushright}

\textbf{Proof of Main Theorem 2} Assuming the degenerate pair $(\mathcal{S}_0,\mathcal{S})$ is not of type $(Z^{n-1}_{\textnormal{max}},Q^{2n-1}),\ \ n\geq 2$ or $(Z^1_{\textnormal{max}},G^{\Rmnum{3}}(n,n)),\ \ n\geq 2$ or $(Q^n,Q^m),\ \ n\equiv m$(mod 1) and following Lemma \ref{Lemma 4.15'}, it is the root vector $E_{\gamma+\theta_{n+1}}$ that causes degeneracy, i.e., $\sigma'(E_{\gamma+\theta_{n+1}},T_{\alpha}(\widetilde{\mathscr{C}}_o(\mathcal{S}_0)))=0$.  Let $(U,\{\omega_j\}^m_{j=1})$ be Harish-Chandra coordinate dual to positive root vectors $E^+$ with the first $n$ coordinate being dual to positive root vectors $E^+_0$. Without loss of generality, let $\omega_{n+1}$ be the coordinate dual to $E_{\gamma+\theta_{n+1}}$. Now consider a germ of submanifold $S\subset\mathcal{S}$ defined on $U$ as image of a mapping as follows: 

\begin{align*}
f:\mathbb{C}^n & \rightarrow\mathbb{C}^n\times\mathbb{C}^{m-n},\\
(z_1,...,z_n) & \rightarrow (z_1,...,z_n,z^2_n,0,...,0),
\end{align*} The tangent vector is of the form $\alpha'+2z_nE_{\gamma+\theta_{n+1}}$. Now we need to show that $S\subset\mathcal{S}$ inherits a sub-VMRT structure modeled on $\mathcal{S}_0$. 

On one hand, if $\alpha'$ is also a VMRT, i.e., $\alpha'\in\widetilde{\mathscr{C}}_o(\mathcal{S}_0)$, then $\alpha'+2z_nE_{\gamma+\theta_{n+1}}\in\widetilde{\mathscr{C}}_o(\mathcal{S})$ by Lemma \ref{Lemma 4.14}.

On the other hand, if $\alpha'+2z_nE_{\gamma+\theta_{n+1}}\in\widetilde{\mathscr{C}}_o(\mathcal{S})$, then $\exists\xi\in T_{\alpha}(\widetilde{\mathscr{C}}_o(\mathcal{S}_0))$ s.t. $\zeta(\xi+2z_nE_{\gamma+\theta_{n+1}})=\alpha'+2z_nE_{\gamma+\theta_{n+1}}$. Since $\sigma'(E_{\gamma+\theta_{n+1}},T_{\alpha}(\widetilde{\mathscr{C}}_o(\mathcal{S}_0)))=\sigma'(E_{\gamma+\theta_{n+1}},E_{\gamma+\theta_{n+1}})=0$, we have $\alpha+(\xi+2z_nE_{\gamma+\theta_{n+1}})+\sigma'(\xi,\xi)=\alpha+(\xi+2z_nE_{\gamma+\theta_{n+1}})+\sigma'(\xi+2z_nE_{\gamma+\theta_{n+1}},\xi+2z_nE_{\gamma+\theta_{n+1}})=\zeta(\xi+2z_nE_{\gamma+\theta_{n+1}})$$=\alpha'+2z_nE_{\gamma+\theta_{n+1}}\Rightarrow\alpha'=\zeta(\xi)\in\widetilde{\mathscr{C}}_o(\mathcal{S}_0)$.

The remaining task to complete the proof of Main Theorem 2 is to take into account the pairs $(Z^{n-1}_{\textnormal{max}},Q^{2n-1}),\ \ n\geq 2$ and $(Z^1_{\textnormal{max}},G^{\Rmnum{3}}(n,n)),\ \ n\geq 2$ and $(Q^n,Q^m),\ \ n\equiv m$(mod 1). For the first two, Hong and Choe (2004) have established example in \cite{[HoC 04]} which shows the non-rigidity. For third one, any $n$-dimensional germ of complex submanifold $S\subset Q^m$ inherits a sub-VMRT structure modeled on $Q^n$, so it is non-rigid. Thus we have completed our proof.

\begin{flushright}
QED
\end{flushright}

\end{document}